\documentclass[12pt]{amsart}
\usepackage{amssymb, amsmath,amsthm}
\usepackage{url}
\newtheorem{thm}{Theorem}[section]
\newtheorem{prop}[thm]{Proposition}
\newtheorem{conj}[thm]{Conjecture}
\newtheorem{cor}[thm]{Corollary}
\newtheorem{lem}[thm]{Lemma}
\theoremstyle{definition}
\newtheorem{rem}[thm]{Remark}
\newtheorem{def1}[thm]{Definition}

\newcommand{\ra}{\rightarrow}
\newcommand{\bk}{\backslash}
\newcommand{\mc}{\mathcal}
\newcommand{\mf}{\mathfrak}
\newcommand{\mb}{\mathbb}
\newcommand{\sg}{\sigma}

\newcommand{\e}{\epsilon}

\newcommand{\mbf}{\boldsymbol}
\renewcommand{\bar}{\overline}
\begin{document}
\title[Bivariate Erd\H{o}s-Kac]{On the Bivariate Erd\H{o}s-Kac Theorem and Correlations of the M\"{o}bius Function}
\author{Alexander P. Mangerel}
\address{Department of Mathematics\\ University of Toronto\\
Toronto, Ontario, Canada}
\email{sacha.mangerel@mail.utoronto.ca}
\maketitle
\begin{abstract}
%Given a positive integer $n$ let $\omega(n)$ denote the number of distinct prime factors of $n$, and let $a$ be fixed positive integer. We develop a bivariate probabilistic model to study the joint distribution of the deterministic vectors $(\omega(n),\omega(n+a))$ with $n \leq x$ as $x \ra \infty$, where $n$ and $n+a$ are both restricted to belong to a subset of $\mb{N}$ with suitable properties. In so doing, we establish a \emph{quantitative} version of a bivariate analogue of the Erd\H{o}s-Kac theorem on proper subsets of $\mb{N}$. \\
%%This is related to the behaviour of the exponential sums $\sum_{n \leq x} \mu^2(n(n+1))e^{i(u\omega(n) + v \omega(n+1))}$, for $u, v \in \mb{R}$. Using this relationship, we show, in an effective way, that these exponential sums have mean value zero in a range of values of $u$ and $v$. \\
%We give two applications of this result. The first is a partial result in the direction of a conjecture of Chowla on binary correlations of the M\"{o}bius function. 
Let $2 \leq y \leq x$ such that $\beta := \frac{\log x}{\log y} \ra \infty$. Let $\omega_y(n)$ denote the number of distinct prime factors $p$ of $n$ such that $p \leq y$, and let $\mu_y(n) := \mu^2(n)(-1)^{\omega_y(n)}$, where $\mu$ is the M\"{o}bius function. We prove that if $\beta$ is not too large (in terms of $x$) then for each fixed $a \in \mb{N}$,
\begin{equation*}
\sum_{n \leq x} \mu_y(n)\mu_y(n+a) \ll x\left(\frac{1}{\log_2 y} + e^{-\frac{1}{21}\beta \log \beta}\right).
\end{equation*}
This can be seen as a partial result towards the binary Chowla conjecture. Our main input is a \emph{quantitative} bivariate analogue of the Erd\H{o}s-Kac theorem regarding the distribution of the pairs $(\omega(n),\omega(n+a))$, where $n$ and $n+a$ both belong to any subset of the positive integers with suitable sieving properties; moreover, we show that the set of squarefree integers is an example of such a set. We end with a further application of this probabilistic result related to a problem of Erd\H{o}s and Mirsky on the number of integers $n \leq x$ such that $\tau(n) = \tau(n+1)$.
%
%We show that, if $\mu_y(n) := \mu^2(n) (-1)^{\omega_y(n)}$, where $\omega_y(n)$ is the number of distinct primes $p|n$ with $p \leq y$, then as long as $\frac{\log x}{\log y} \ra \infty$ as $x \ra \infty$,
%\begin{equation*}
%\sum_{n \leq x} \mu_y(n)\mu_y(n+a) = o(x).
%\end{equation*}
%The second application is a partial result in the direction of a conjecture of Erd\H{o}s, Pomerance and Sark\H{o}zy on the order of magnitude of the number of $n \leq x$ such that $\tau(n) = \tau(n+1)$, where $\tau$ is the divisor function. We establish that if $\tau_y(n)$ is the number of $y$-smooth divisors of $n$ with $y \leq x^{\frac{1}{\log_3 x}}$ then for all $|j| \leq \left(1-\e\right)\sqrt{\log_3 x}$,
%\begin{equation*}
%|\{n \leq x : \tau_y(n) = 2^j\tau_y(n+1)\}| \gg e^{-\frac{j^2}{4}}\frac{x}{\sqrt{\log_2 x}},
%\end{equation*}
%which is the conjectured order of magnitude of the lower bound when $j = 0$ and $y = x$.
\end{abstract}
\section{Introduction}
\subsection{On the Binary Chowla Conjecture}
Let $\mu$ denote the M\"{o}bius function. It is well-known that the Prime Number Theorem is equivalent to the statement that $\sum_{n \leq x} \mu(n) = o(x)$. Thus, $\mu$ exhibits a lot of cancellation, and we in fact expect (according to the density hypothesis) that for each $\e > 0$ and $x$ sufficiently large, $\mu$ has a sign change in any interval of the form $(x,x+x^{\e}]$. We also expect that the sign changes are random, and do not exhibit conspiratorial tendencies, such as $\mu(n)$ and $\mu(n+a)$ frequently changing sign simultaneously for a given, fixed $a \in \mb{N}$. One of the first enunciations of this latter principle is due to Chowla \cite{Cho}.
\begin{conj}[Chowla] \label{CONJCHOWLA}
Let $k \geq 2$ and let $0 \leq a_1 < a_2 < \cdots < a_k$ be integers. Then
\begin{equation} \label{FULLCHOW}
\sum_{n \leq x} \mu(n+a_1)\cdots \mu(n+a_k) = o(x).
\end{equation}
\end{conj}
In the binary case, i.e., $k = 2$, Chowla's conjecture is the statement that for any fixed $a \in \mb{N}$,
\begin{equation} \label{BINCHOW}
\sum_{n \leq x} \mu(n)\mu(n+a) = o(x).
\end{equation}
This is currently not known for \emph{any} $a \in \mb{N}$. In fact, even to show that the left side of \eqref{BINCHOW} has absolute value at most $cx$, where $c < 1$, was an intractable problem until very recently, when Matom\"{a}ki and Radziwi\l \l $ \ $ were able to prove this as a consequence of their work on short averages of multiplicative functions (see Corollary 2 of \cite{MaR} for the corresponding result with the Liouville function $\lambda$ in place of $\mu$). \\
There has been some remarkable recent progress on Conjecture \ref{CONJCHOWLA} itself. We mention the two following notable examples. Tao \cite{Tao} proved a logarithmically averaged version of \eqref{BINCHOW}, that is
\begin{equation} \label{LOGWEIGHT}
\sum_{n \leq x} \frac{\mu(n)\mu(n+a)}{n} = o(\log x).
\end{equation}
%whenever $g_1,g_2$ are non-pretentious, in the sense of Granville and Soundararajan \cite{Tao}. Since $\mu$ is non-pretentious, Tao's theorem implies a log-weighted version of Chowla's conjecture. 
Unfortunately, \eqref{LOGWEIGHT} is a strictly weaker estimate than \eqref{BINCHOW}. In a different direction, Matom\"{a}ki, Radziwi\l \l $\ $and Tao \cite{MRT} proved that if one averages over all shift vectors $\mbf{a} \in [0,H]^k$, where $H \ra \infty$ as $x \ra \infty$ then \eqref{BINCHOW} holds for almost every such $\mbf{a}$.  \\
We shall prove the following two partial results in the direction of \eqref{BINCHOW}, which provide further motivation for the binary case of Conjecture \ref{CONJCHOWLA} for any \emph{fixed} shift. \\
Set $B(x) := \frac{\sqrt{\log_2 x}}{\log_3 x \log_4^2 x}$, where $\log_k x$ is the $k$th iterated logarithm and $x$ is sufficiently large.
%From a heuristic standpoint this is quite reasonable. Indeed, $\phi_x$ is the characteristic function for the law of the distribution $F_x$. LeVeque \cite{LEV} proved that (if we remove the condition $\mu^2(n(n+1)) = 1$ in the definition of $F_x$) a joint Erd\H{o}s-Kac theorem holds, i.e., the $L^{\infty}$ distance between $F_x$ and the distribution function of a bivariate Gaussian tends to 0 as $x \ra \infty$. It is well-known that the convergence in distribution of a sequence of random vectors to a limit random vector is equivalent to the \emph{pointwise} convergence of the corresponding characteristic functions of these random vectors. Thus, modulo the difference in the arithmetic conditions in the definition of $F_x$ and $\phi_x$, we know roughly that $|\phi_x(u,v)-\chi(u,v)|$ vanishes as $x \ra \infty$ for \emph{fixed} $u,v$; the assertion that this continues to be true for $u,v$ a function of $x$ is a more subtle and difficult one. While we do not demonstrate that the convergence is uniform, we prove in the sequel an effective bound, in terms of $u$ and $v$, for this distance, a corollary of which is the following.
\begin{thm}\label{CHOW}
a) Let $2 \leq y \leq x$ and let $\beta := \frac{\log x}{\log y}$. Suppose that $\beta \leq e^{B(x)}$. For $n \in \mb{N}$ let $\omega_y(n)$ denote the number of primes $p$ dividing $n$ with $p \leq y$, and let $\mu_y(n) := \mu^2(n)(-1)^{\omega_y(n)}$.  Then for each fixed $a \in \mb{N}$,
\begin{equation}\label{ONLYZ}
\sum_{n \leq x} \mu_y(n)\mu_y(n+a) \ll x\left(\frac{1}{\log_2 y} + e^{-\frac{1}{21} \beta \log \beta}\right).
%\begin{cases} \frac{x \log_4 x}{\sqrt{\log_2 x \log_3 x}} &\text{ if $\beta \geq \log_3 x$} \\ .
\end{equation}
b) Let $\mu(n;u) := \mu^2(n)e^{iu\omega(n)}$. If $|u|,|v| \leq 2\pi$ and $w := \max\{|u|,|v|\}$, we have
\begin{equation} \label{SMALLU}
%\sum_{n \leq x} \mu^2(n)\mu^2(n+1) e^{iu\omega_x(n)} e^{iv\omega_x(n+1)} 
\sum_{n \leq x} \mu(n;u)\mu(n+a;v) \ll x\left(w\log(1/w) + (\log x)^{-\frac{1}{3}(u^2+v^2)}\right).
\end{equation}
\end{thm}
Note that \eqref{ONLYZ} states that if we only account for those prime factors $p$ of integers $n \leq x$ such that $p \leq x^{o(1)}$ then the corresponding correlation sums of $\mu_y$ are small. In other words, those sign changes of $\mu(n)$ that are caused by ''small primes'' do not appear to correlate with those of $\mu(n+a)$. \eqref{SMALLU}, while not directly related to \eqref{BINCHOW}, roughly shows that if we replace $\mu(n) = \mu^2(n)e^{i\pi \omega(n)}$ by $\mu^2(n)e^{iu\omega(n)}$ for \emph{any} $u = o(1)$ that is not too small as $x \ra \infty$ then the corresponding correlation sums are also small. 
\begin{rem}
Results like \eqref{SMALLU} appear to be completely new, and provide a collection of examples in the direction of a general conjecture on the size of correlations of non-pretentious multiplicative functions, originally due to Elliott (see, for example, Conjecture 1.5 in \cite{MRT}). \\
Results like \eqref{ONLYZ} do exist in the literature, but in weaker forms. To the author's knowledge, the first result of this form is Theorem 5 in \cite{CFMRS}, where a correlation estimate of the type of \eqref{ONLYZ} is established for a truncated version of the Liouville function, i.e., $\lambda_y(n) := (-1)^{\Omega_y(n)}$, where $\Omega_y(n)$ is the number of prime factors $p$ of $n$ with $p \leq y$, counted with multiplicity. However, the parameter $y$ is restricted in such a way that $y \leq (\log x)^{2+o(1)}$. \\ 
Subsequently, Daboussi and Sark\H{o}zy \cite{DaS} established a more general theorem, applicable to the truncation of \emph{any} 1-bounded multiplicative function, in which one may take any $y \leq x^{o(1)}$. As a particular case, they show that
\begin{equation*}
\sum_{n\leq x} \mu_y(n)\mu_y(n+1) \ll x\left(\frac{1}{(\log y)^{9}} + e^{-\frac{\beta}{8}}\right).
\end{equation*}
Their result is not explicitly stated for shifts $a$ other than 1, though this restriction appears to be merely technical. We emphasize, though, that our probabilistic view of the problem motivates the methods that we employ in this paper, and these are substantially different from those used in \cite{DaS}. Among other things, in \cite{DaS} Brun's sieve is used, while in the present paper we appeal to the Rosser-Iwaniec sieve instead (see Lemma \ref{ROSSER} below). This difference accounts for the appearance of the additional factor of $\log \beta$ in the second term in \eqref{ONLYZ} above. In particular, when $y \geq x^{\frac{1}{3\log_3 x}}$ our estimate is superior to theirs. We believe that this improvement is worthwhile, considering that 
%
%We emphasize here that, aside from their result being different than ours and that our method is completely different from theirs, they can only take $y \leq (\log x)^{2+o(1)}$. 
%This is less satisfying than \eqref{ONLYZ}, in the sense that 
the main interest in a result like \eqref{ONLYZ} is in its effectiveness when $y$ is as close to $x$ as possible in support of Conjecture \ref{CONJCHOWLA}. \\
% and Theorem \ref{CHOW} provides \emph{any} $y \leq x^{o(1)}$ that is sufficiently large. 
We also note that a ternary extension of the result in \cite{DaS} was worked out by Ganguli in his thesis \cite{Gan}. This extension essentially follows the same method of proof as that in \cite{DaS}.
\end{rem}
As our results suggest, the main focus of our arguments is on the behaviour of the pairs $(\omega(n),\omega(n+a))$, where $n$ and $n+a$ are squarefree. Indeed, observe that we can express \eqref{BINCHOW} equivalently in the form
\begin{equation*}
\sum_{n \leq x} \mu^2(n)\mu^2(n+a) e^{iu \omega(n)}e^{iv \omega(n+a)} = o(x),
\end{equation*}
where $u = v = \pi$. A moment's reflection (assuming one has probabilistic inclinations) suggests that the function of $u$ and $v$ resembles some variant of a characteristic function for a random vector. That is, suppose $\mbf{X} := (X_1,X_2)$ is a random vector on a fixed probability space $(\Omega,\mc{B},\mb{P})$, and let $\sigma_{\mbf{X}} := \mb{P} \circ \mbf{X}^{-1}$ denote its law. We recall that the \emph{characteristic function} of $\mbf{X}$ is the Fourier transform of $\sigma_{\mbf{X}}$, i.e., 
\begin{equation*}
\phi_{\mbf{X}}(\mbf{t}) := \mb{E}[e^{i \mbf{X} \cdot \mbf{t}}] = \int_{\mb{R}^2} e^{i \mbf{u} \cdot \mbf{t}} d\sigma_{\mbf{X}}(\mbf{u}).
\end{equation*}
Let $S_a$ denote the set of integers $n$ such that $n$ and $n+a$ are both squarefree. Choosing the finite probability space given by $[1,x] \cap S_a$ with its power set and normalized counting measure, the characteristic function of the random vector $n \mapsto (\omega(n),\omega(n+a))$ is
\begin{equation*}
\tilde{\phi}_x(\mbf{t}) := E^{\ast}(x;a)^{-1}\sum_{n \leq x} \mu^2(n)\mu^2(n+a) e^{i(\omega(n),\omega(n+a))\cdot \mbf{t}},
\end{equation*}
where $E^{\ast}(x;a) := \sum_{n \leq x} \mu^2(n)\mu^2(n+a)$. In particular, \eqref{BINCHOW} is stating that $\tilde{\phi}_x(\pi,\pi) = o(1)$, as $x \ra \infty$. In the next subsection, we shall discuss an analogue of this estimate that \emph{is} achievable, by way of a quantitative, bivariate generalization of the Erd\H{o}s-Kac theorem. This is the crucial input into the proof of Theorem \ref{CHOW}.
\begin{rem}
In principle, our arguments extend to 
%more general correlations $\mu_z(n)\mu_z(an+b)$, with $a,b \in \mb{N}$. With more effort, our arguments should also be applicable to 
$k$-ary correlations of $\mu_y$ and $\mu(\cdot;u)$, as well as to certain classes of unimodular multiplicative functions. We postpone such extensions to a separate paper.
\end{rem}
\begin{rem}
In a sense, the approach we use as support for Chowla's conjecture is misguided. In light of an argument of Tao \cite{BLOG}, it is known that the estimate
\begin{equation}
\sum_{n \leq x} \mu(n)\mu(an+2) \ll_{\e} \frac{x}{(\log x)^{2+\e}}, \label{STRONGCHOW}
\end{equation}
for each $a \in \mb{N}$, would be sufficient to prove the Twin Prime conjecture. The notorious Parity Problem in Sieve Theory prevents the Twin Prime conjecture from being tractable via Sieve methods alone; hence, it would seem that an attempt at proving Chowla's conjecture must also invoke additional parity barrier-breaking arguments. Our approach, as outlined in Section 2, involves generalizing a probabilistic framework of Kubilius in order to analyze joint distributions, and makes use of a composition of Rosser-Iwaniec sieves (see Lemma \ref{ROSSER}). Thus, our method is heuristically insufficient to provide a full proof of Chowla's conjecture. In spite of this, we are content to provide here a basic framework for correlation problems upon which further investigations may build.
\end{rem}
\subsection{A Quantitative Bivariate Erd\H{o}s-Kac Theorem}
%A major influence in the development of Probabilistic Number Theory is the Erd\H{o}s-Kac theorem (see Section III.4 in \cite{Ten2}). 
Given $n \leq x$, let $$\tilde{\omega}_x(n) := (\omega(n)-\log_2 x)/\sqrt{\log_2 x},$$ where $\omega(n)$ is the number of distinct prime factors of $n$. Furthermore, for $z \in \mb{R}$ let 
\begin{equation*}
H_x(z) := x^{-1}|\{n \leq x : \tilde{\omega}_x(n) \leq z\}|.
\end{equation*}
$H_x$ is the distribution function for $\tilde{\omega}_x(n)$, a centred and normalized version of $\omega$. It is well-known that if we define a collection $\{\theta_p\}_p$ of indicator functions indexed by primes such that $\theta_p(n) = 1$ if $p|n$ and $\theta_p(n) = 0$ otherwise then $\omega(n) = \sum_{p \leq x} \theta_p(n)$ for every $n \leq x$, and the values of $\theta_p$ and $\theta_q$ are \emph{asymptotically} independent as $x \ra \infty$. In spite of this, $\omega$ still behaves like a sum of \emph{genuinely} independent random variables in the sense that it obeys a Central Limit Theorem. Indeed, this conclusion is furnished by the Erd\H{o}s-Kac theorem, which states that as $x \ra \infty$, $H_x \ra \Phi$ almost everywhere, where 
\begin{equation*}
\Phi(z) := \frac{1}{\sqrt{2\pi}}\int_{-\infty}^z e^{-\frac{1}{2}t^2} dt  
\end{equation*}
is the distribution function of a Gaussian random variable.  The rate of convergence of $H_x$ to the normal distribution has also been studied in this connection, and a best-possible estimate for the $L^{\infty}$ distance between $H_x$ and $\Phi$ was conjectured by LeVeque \cite{LEV} and proven by R\'{e}nyi and Tur\'{a}n \cite{ReTu} to be $O\left(1/\sqrt{\log_2 x}\right)$. This result echoes the best-possible rate in the Lindeberg Central Limit theorem, given by the Berry-Ess\'{e}en Theorem (see Chapter XVI of \cite{FEL}). \\
One expects that since $n$ and $n+a$ share only finitely many common prime factors, the values of $\tilde{\omega}_x(n)$ and $\tilde{\omega}_x(n+a)$ are also asymptotically independent as $x \ra \infty$. It is therefore reasonable to guess that a \emph{bivariate} analogue of the Erd\H{o}s-Kac theorem should hold, in which the distribution function $H_x$ is replaced by the two-dimensional analogue
\begin{equation*}
H'_x(z,z') := x^{-1} \left|\left\{n \leq x: \tilde{\omega}_x(n) \leq z, \tilde{\omega}_x(n+a) \leq z'\right\}\right|,
\end{equation*}
and the Gaussian limit distribution is replaced by the uncorrelated bivariate Gaussian distribution $$\Phi_{(2)}(z,z') := \Phi(z)\Phi(z').$$ This was in fact proven by LeVeque \cite{LEV}, though his result does not yield an effective rate of decay for $\|H'_x - \Phi_{(2)}\|_{L^{\infty}(\mb{R}^2)}$.\\
\begin{rem}
Actually, more general theorems than that of LeVeque regarding the existence and characterization of joint distributions of additive functions and their shifts appear in the literature; for a synthesis of these results, see Chapter VII of \cite{KUB}. However, the methods employed to prove these results, e.g., the Cr\'{a}mer-Wald trick, are typically not useful for producing quantitative results. In a related though distinct vein, see Chapter V in \cite{KUB} for quantitative results related to the distribution of sums $\sum_{1 \leq j \leq k} f_j(n+a_j)$, where each $f_j$ is a real-valued, strongly multiplicative function.
\end{rem}
In the sequel, we will concern ourselves with questions surrounding a restricted analogue of LeVeque's theorem. Specifically, let $a \in \mb{N}$ be fixed, let $E^{\ast}(x;a)$ and $\tilde{\omega}_x$ be as above and for $z,z' \in \mb{R}$ let
\begin{equation*}
F_{x,a}(z,z') := E^{\ast}(x;a)^{-1} \left|\{n \leq x : \mu^2(n) = \mu^2(n+a) = 1, \tilde{\omega}_x(n) \leq z, \tilde{\omega}_x(n+a) \leq z'\}\right|.
\end{equation*}
We shall prove the following.
\begin{thm}\label{MAINEK}
Let $a \in \mb{N}$. Then as $x \ra \infty$,
\begin{equation*}
\|F_{x,a} - \Phi_{(2)}\|_{L^{\infty}(\mb{R}^2)} \ll (\log_3 x)(\log_2 x)^{-\frac{1}{4}}.
\end{equation*}
\end{thm}
We will actually prove a rather more general result, Theorem \ref{EKTHM}, which gives an $L^{\infty}$ distance estimate as in Theorem \ref{MAINEK} for a distribution function associated to pairs $(\omega(n),\omega(n+a))$, where $n$ and $n+a$ belong to \emph{any} subset of $\mb{N}$ with suitable sieve properties. Theorem \ref{MAINEK} is a corollary of this (in this connection, see Proposition \ref{QUADRATFREI}).
\begin{rem}
As in the univariate case, we expect that the best error term here is $O\left(1/\sqrt{\log_2 x}\right)$ (to see that this is best possible, see Lemma \ref{CONTRA}). As mentioned, R\'{e}nyi and Tur\'{a}n \cite{ReTu} obtained this error term in the classical Erd\H{o}s-Kac theorem using analytic methods from multiplicative number theory. Indeed, their proof depends crucially on the fact that the characteristic function of the distribution they consider is
\begin{equation*}
u \mapsto x^{-1} \sum_{n \leq x} e^{iu \omega(n)},
\end{equation*}
the mean value of the \emph{multiplicative} function $n \mapsto e^{iu\omega(n)}$ for $n \leq x$. This can be estimated quite precisely for each $u \in [0,2\pi]$ by the Selberg-Delange method (see Chapter II.5 of \cite{Ten2}). Such techniques are not at our disposal, however, and our methods are not sufficiently powerful to yield the best-possible estimate.
\end{rem}
\begin{rem}
As an application of the more general Theorem \ref{EKTHM} below, one can deduce an asymptotic formula for the number of integers $n \leq x$ (\emph{without} the squarefree restriction) such that $\omega(n) = k_1$ and $\omega(n+a) = k_2$, where each $k_j$ is sufficiently close to $\log_2 x$, and $a \in \mb{N}$ is fixed. We note in this connection that Goudout (see Th\'{e}or\`{e}me 3 in \cite{Gou}) has recently proved an upper bound of the correct order of magnitude for the number of such integers that is uniform over all $1 \leq k_1,k_2 \leq R \log_2 x$, given any fixed $R > 0$.  
\end{rem}
It will be convenient for us to use a slight modification of the function $\tilde{\phi}_x$ defined in the previous subsection. For $u,v \in \mb{R}$ let
\begin{equation*}
\phi_{x,a}(u,v) := E^{\ast}(x;a)^{-1}\sum_{n \leq x} \mu^2(n)\mu^2(n+a)e^{i(u\tilde{\omega}_{x}(n) + v\tilde{\omega}_{x}(n+a))},
\end{equation*}
In terms of $\phi_{x,a}$, \eqref{BINCHOW} is equivalent to the statement that $\phi_{x,a}(u,v) = o(1)$, where $u = v = \pi \sqrt{\log_2 x}$. Note that for this choice of $u$ and $v$, $\chi(u,v) = o(1)$, so it would suffice to know that $|\phi_{x,a}(u,v) - \chi(u,v)| = o(1)$ to prove \eqref{BINCHOW}.
%and let $\chi(u,v) := e^{-\frac{1}{2}(u^2+v^2)}$. For convenience, we will write $\phi_x$ when $S$ and $a$ are understood. It is easy to see that $\phi_x(u,v)$ and $\chi(u,v)$ are the respective Fourier transforms of the Stieltjes measures $dF_{x,S,a}(z,z')$ and $d\Phi_{(2)}(z,z')$; $\phi_x$ and $\chi$ are known, respectively, as the \emph{characteristic functions} of $F_{x,S,a}$ and $\Phi_{(2)}$. 
By a well-known theorem of L\'{e}vy, the convergence in distribution $\|F_{x,a}-\Phi_{(2)}\|_{L^{\infty}(\mb{R}^2)} \ra 0$ is equivalent to the \emph{pointwise} convergence of the corresponding characteristic functions of these distributions, i.e., $|\phi_{x,a}(u,v)-\chi(u,v)| \ra 0$ as $x \ra \infty$ for each $u,v \in \mb{R}$. The assertion that this continues to be true when $u,v$ are chosen to \emph{grow as a function of $x$} is a more subtle and difficult one. While we do not demonstrate that the convergence is uniform in general, we develop a method to prove an effective bound, in terms of $u$ and $v$, for the distance between $\phi_{x,a}(u,v)$ and $\chi(u,v)$ in a range of $u$ and $v$ depending on $x$ (that unfortunately does not include $u = v= \pi\sqrt{\log_2 x}$). We shall use this bound in conjunction with a bivariate analogue of a smoothing lemma of Ess\'{e}en (see Lemma \ref{ESSEEN}) to prove Theorem \ref{MAINEK}.
%, up to an extra $\log_3 x$ factor; the exponent 2 on the triple log factor can be seen to reflect the fact that the joint law of two random variables, rather than a law of a single variable, is in question here.\\
%Our interest is not purely in whether or not the bivariate Central Limit phenomenon persists in proper subsets of $\mb{N}$; rather, we are also motivated by how this phenomenon relates to the two applications that we list in the subsections below, and our methods happen to lend themselves to the proof of Theorem \ref{EKTHM} as well. 
%One consequence of our analysis is the following.
%The method we present here is limited by techniques from Probability theory, specifically a theorem of Roos \cite{ROOS} on multivariate Poisson approximation (see the end of Section 3 for details). It is possible that a more careful analysis of the rate of approximation in our specific case might yield an estimate like \eqref{STRONGCHOW}.
\subsection{On a Problem of Erd\H{o}s and Mirsky}
For a positive-valued multiplicative function $f$ let
\begin{equation*}
S_f(x) := |\{n \leq x: f(n) = f(n+1)\}|.
\end{equation*}
It is generally a difficult problem to even determine whether $S_f(x)$ tends to infinity with $x$. Much of the literature on problems of this type relate to the case where $f = \tau$, the divisor function, and this shall be our focus as well. We thus henceforth write $S(x)$ to mean $S_{\tau}(x)$. \\
Erd\H{o}s and Mirsky \cite{EM} famously conjectured that $\lim_{x \ra \infty} S(x) = \infty$, i.e., that there are infinitely many integers $n$ such that $\tau(n) = \tau(n+1)$. Based on ideas of C. Spiro, Heath-Brown proved this conjecture in the affirmative, giving the lower bound $S(x) \gg x(\log x)^{-7}$. Erd\H{o}s, Pomerance and Sark\H{o}zy have made the following conjecture regarding the order of magnitude of $S(x)$.
\begin{conj}[\cite{EPS}] \label{CONJEPS}
We have $S(x) \asymp x(\log_2 x)^{-\frac{1}{2}}$.
\end{conj}
The latter three authors proved that the upper bound in Conjecture \ref{CONJEPS} holds, and cite an heuristic argument due to Bateman and Spiro as further motivation for this conjecture. In the opposite direction, Hildebrand \cite{Hil} has shown that $S(x) \gg x(\log_2 x)^{-3}$. See the beginning of Section 6 for our version of the Bateman-Spiro heuristic. \\
%, based on, among other things, an upper bound argument due to Erd\H{o}s, Pomerance and Sark\H{o}zy \cite{EPS}. In the opposite direction, 
%, although the Bateman-Spiro heuristic is not explicitly described in their paper. \\
Implementing techniques related to those used to prove Theorem \ref{CHOW} and Theorem \ref{EKTHM}, we can prove a partial result in this direction. 
% of the conjecture of Erd\H{o}s, Pomerance and Sark\H{o}zy. 
%
%and will in fact prove the validity of the heuristic, and thus the conjecture of Erd\H{o}s-Pomerance-Sark\H{o}zy.  
Let $A(x) := 21\frac{\log_3 x}{\log_4 x}$.
\begin{thm}\label{MAINEM}
Let $2 \leq y \leq x$ such that if $y = x^{\frac{1}{\beta}}$ then $A(x) \leq \beta \leq e^{B(x)}$. Let $\tau_y(n)$ be the number of divisors $d |n$ such that if $p|d$ then $p \leq y$. 
%Then for each $|j| \leq \left(1-\e\right)\sqrt{\log_3 x}$, we have
\begin{equation*}
%|\{n \leq x : \tau_y(n) = 2^j\tau_y(n+1)\}| \gg  e^{-\frac{j^2}{4}}\frac{x}{\sqrt{\log_2 x}}.
|\{n \leq x : \tau_y(n) = \tau_y(n+1)\}| \gg  \frac{x}{\sqrt{\log_2 x}}.
\end{equation*}
\end{thm}
In Section 6 we actually prove a more general result on the number of $n \leq x$ such that $\tau(n) = 2^j \tau(n+1)$ in a range of $j$ depending on $x$ that includes $j = 0$. See Theorem \ref{EM} for a precise statement. \\
%Note that the set of $n \leq x$ such that $\tau(n) = \tau(n+1)$ distributes among sets in which $\tau_y(n) = 2^j \tau_y(n+1)$, for some $|j| \leq \frac{\log x}{\log y}$. However, most integers have about $\log\left(\frac{\log x}{\log y}\right)$ such factors, and we thus expect that $n$ and $n+a$ should have nearly the same number of prime factors of size larger than $y$. Thus, the restriction on $j$ here is not expected to be problematic when $y$ is close to the upper limit of its range. If we had a better understanding of which $j$ occur most often such that $\tau(n) = \tau(n+1)$ \emph{and} $\tau_y(n) = 2^j\tau_y(n+1)$ then we would have a chance at actually proving Conjecture \ref{CONJEPS}.\\ 
Our methods also extend to the more general question of estimating from below the number of integers $n \leq x$ with $\tau_y(n) = \tau_y(n+a)$ for $a \in \mb{N}$, and with more effort, to questions such as determining how often $\tau_y(n+a_1) = \cdots = \tau_y(n+a_k)$, for $0 \leq a_1 < \cdots < a_k$ and $k \geq 2$. 
\subsection{Notation and Conventions}
Throughout this paper, $a,b,k,l,m,d$ will always stand for positive integers, and $p,q$ will always denote primes. We write $\log_k x$ to mean $\log(\log_{k-1} x)$, with $\log_1 x = \log x$, and $x$ is always assumed to be sufficiently large so that these quantities are \emph{positive}. We will frequently write $\mathfrak{L}_y$ to mean $\sqrt{\log_2 y}$ and $\mf{L} := \mf{L}_x$. We will also denote by $\sg$ the shift map $\sg(n) := n+1$ defined on $\mb{N}$.\\
We shall employ the usual conventions of Probability Theory. We denote by $\mb{P}$ a fixed probability measure on a measurable space $(\Omega,\mc{B})$ on which a given random vector $X: \Omega \ra \mb{R}^n$ is defined. If $S \subseteq \mb{R}^n$ is Borel measurable and $X^{-1}(S) =: A \in \mc{B}$ then we will write $\mb{P}(X \in S)$ in place of $\mb{P}(A)$. We write $\mc{L}(X)$ to denote the measure $\mb{P} \circ X^{-1}$, which we call the \emph{law} of $X$. Obviously, the law of any random vector is a probability measure on $\mb{R}^n$ equipped with its Borel $\sg$-algebra. We let $\mb{E}X$ denote the expectation of $X$, i.e., 
\begin{equation*}
\mb{E}X := \int_{\Omega}X d\mb{P} = \int_{\mb{R}^n} \mbf{t}d\mc{L}(X)(\mbf{t}). 
\end{equation*}
Given two probability measures $\nu_1$ and $\nu_2$ on the measurable space $(\Omega,\mc{B})$ we define the \emph{total variation distance} between them by
\begin{equation*}
d_{TV}(\nu_1,\nu_2) := \sup_{A \in \mc{B}} |\nu_1(A)-\nu_2(A)|. 
\end{equation*}
\section{A Bivariate Kubilius Model and Multivariate Poisson Approximation}
%In this connection we introduce the following definition.
In this section we give a general framework to study the distribution of the vector $(\omega(n),\omega(n+a))$, where $n$ and $n+a$ are confined to subsets of $\mb{N}$ with suitable properties. To be precise, we introduce the following definition. 
\begin{def1} \label{SIFTABLE}
Let $a \in \mb{N}$. Let $S \subseteq \mb{N}$. We say that $S$ is \emph{siftable} with respect to $a$ if: a) we have $$E(x;a) := \sum_{n \leq x \atop (n,a) = 1} 1_S(n) 1_S(n+a) \gg_a x,$$
where $1_S$ is the indicator function for $S$, and:  b) there exist a non-negative multiplicative function $f = f_a$ and a real number $\theta \in (0,1)$ 
%and $A \geq 1$ 
such that: \\
i) $0 \leq f(p) \leq \frac{1}{p}$ for each prime $p$, \\
ii) for each pair of coprime squarefree integers $k_1,k_2$ with $k_1k_2 \leq x$ we have
\begin{equation}\label{SIEVE}
\sum_{n \leq x \atop n \equiv 0 (k_1), n + a\equiv 0 (k_2)} 1_S(n)1_S(n+a) 1_{(n,a) = 1} = f(k_1k_2)\left(E(x;a) + R(x;k_1,k_2)\right),
\end{equation}
and the remainders $R(x;k_1,k_2)$ satisfy
\begin{equation*}
\sum_{k_1k_2 \leq x^{\theta}} |R(x;k_1,k_2)| \ll \frac{x}{(\log x)^3}.
\end{equation*}
Finally, we will say that a siftable set is \emph{regular} if $\sum_{p \leq y} f(p) \gg \log_2 y$ for all $y$ sufficiently large.
\end{def1}
Without loss of generality, we can, and will, assume that $f(p) = 0$ whenever $p | a$.\\
Trivially, $S = \mb{N}$ is siftable and regular. Less trivially:
\begin{prop}\label{QUADRATFREI}
The set of squarefree integers is siftable and regular with respect to each fixed $a \in \mb{N}$.
\end{prop}
We shall prove this fact in Section 4, as it is crucial in our application to both Theorems \ref{CHOW} and \ref{EKTHM}.
%It will be relevant in the applications that we shall state below. 
By a similar argument, the set of $r$-free integers for each $r \geq 3$ is siftable. We leave the verification of this statement to the reader.\\
Given a siftable set $S$ and a fixed $a \in \mb{N}$, let $E(x;a)$ denote the number of $n \leq x$ coprime to $a$ such that both of $n,n+a \in S$, and let $E^{\ast}(x;a)$ denote the number of such $n \leq x$ without the coprimality condition. In this context, put $\lambda(x) := \sum_{p \leq x} f(p)$ and
\begin{equation*}
\tilde{\omega}_{x,S}(n) := \frac{\omega(n)-\lambda(x)}{\sqrt{\lambda(x)}},
\end{equation*}
and set
\begin{align*}
F_{x,S,a}(z,z') &:= E(x;a)^{-1} \left|\left\{n \leq x: n,n+a \in S, (n,a) = 1, \tilde{\omega}_{x,S}(n) \leq z, \tilde{\omega}_{x,S}(n+a) \leq z'\right\}\right|,\\
F_{x,S,a}^{\ast}(z,z') &:= E^{\ast}(x;a)^{-1} \left|\left\{n \leq x: n,n+a \in S, \tilde{\omega}_{x,S}(n) \leq z, \tilde{\omega}_{x,S}(n+a) \leq z'\right\}\right|.
\end{align*}
%where $E(x) := \sum_{n \leq x} \mu^2(n(n+1))$ (it is well-known, see e.g., \cite{HB}, that $E(x) \gg x$). 
%We have restricted ourselves to a count of integers $n$ such that $n$ and $n+1$ are both squarefree. 
For convenience, we will write $F_x$ and $F_x^{\ast}$ in place of $F_{x,S,a}$ and $F_{x,S,a}^{\ast}$, respectively, whenever $S$ and $a$ are clearly defined. We shall establish an estimate for the rate of convergence of the limiting process $F_{x,S,a}(z,z') \ra \Phi(z)\Phi(z')$ as $x \ra \infty$. In the statement below we denote by $\sg$ the shift map $\sg(n) := n+1$, for $n \in \mb{N}$, and we write $\sg^a$ to denote the composition of $\sg$ with itself $a$ times.
\begin{thm}\label{EKTHM}
Let $a \in \mb{N}$ be fixed and let $S$ be a siftable, regular set with respect to $a$. Then $$\|F_{x,S,a}-\Phi_{(2)}\|_{L^{\infty}(\mb{R}^2)} \ll (\log_3 x)(\log_2 x)^{-\frac{1}{4}}.$$ \\
Moreover, if $1_S$ is multiplicative and $S$ is also siftable for each divisor of $a$ then $$\|F_{x,S,a}^{\ast} - \Phi_{(2)}\|_{L^{\infty}(\mb{R}^2)} \ll (\log_3 x)(\log_2 x)^{-\frac{1}{4}}.$$ In particular, the limiting distribution of the vectors $(\tilde{\omega}_{x,S},\tilde{\omega}_{x,S} \circ \sg^a)$ in the set $S \cap (S-a)$ is the uncorrelated bivariate Gaussian distribution. 
\end{thm}
Let $S$ be a siftable set and let $a \in \mb{N}$ be fixed. 
%We shall construct a bivariate random model for pairs of consecutive integers both of which are squarefree.
%\begin{my
%Let $S \subseteq \mb{N}$ be a set with positive natural density, and let $S(x) := |S \cap [1,x]|$. 
%Recall that $$E(x;a) := \sum_{n \leq x \atop (n,a) = 1} 1_S(n)1_S(n+a).$$ By definition, there is a multiplicative function $f : \mb{N} \ra [0,1]$ and a number $\theta > 0$ (both possibly depending on $a$)  such that for $k_1k_2 \leq x^{\theta}$, $(k_1,k_2) = 1$,
%\begin{equation} \label{APPROXFORM}
%\sum_{n \leq x \atop n \equiv 0 (k_1), n + a \equiv 0 (k_2)} 1_S(n)1_S(n+a)1_{(n,a) = 1} = f(k_1k_2)\left(E(x;a) + R(x;k_1,k_2)\right),
%\end{equation}
%and the remainders $|R(x;k_1,k_2)|$ are small on average (see Definition \ref{SIFTABLE}). 
%as in Corollary \ref{ROSSCOR}. 
%Let $f(p) := \frac{1}{p}\left(1-\frac{1}{p}\right)\left(1-\frac{2}{p^2}\right)^{-1}$, and extend $f$ multiplicatively to squarefree integers. Then, for $(k_1,k_2) = 1$, $f$ satisfies the statement
%\begin{equation} \label{SIEVE}
%\sum_{n \leq x \atop n \equiv 0(k_1), n+1\equiv 0(k_2)} \mu^2(n(n+1)) = Xf\left(k_1k_2\right) + R(x;k_1,k_2),
%\end{equation}
%where $X = E(x)$ and $R(x;k_1,k_2)$ is the error in Lemma \ref{SIEVESTUFF} with $d = 1$, $q = k_1$ and $r = k_2$. 
Let $y := x^{\frac{1}{\beta}}$ where $\beta \geq 1$, and put $P := P(y) := \prod_{p \leq y \atop p \nmid a} p$.  Our goal is to approximate the vector $(\omega,\omega \circ \sg^a)$ on $S \cap (S-a)$ by a discrete random vector. In this direction, it is sufficient to  construct a measure against which there is a collection of \emph{independent} random vectors $(\tilde{\theta}_p,\tilde{\theta}_q)$that are good approximations for the deterministic events $\{n\leq x : n,n+a \in S, p|n, q|n+a\}$, for $p,q$ distinct primes dividing $P$. \\
Given coprime $k_1,k_2|P$, define the set
%The algebra of sets on which this measure is defined is generated by sets
\begin{equation*}
E_{k_1,k_2} := \left\{n \leq x : (n,a) = 1, n,n+a \in S, k_1|n, k_2|(n+a), (n,P/k_1) = (n+a,P/k_2) = 1\right\}.
\end{equation*}
%The constructed measure alluded to will be an approximation of the measure $\nu_z$ whose values on the generating set $E_{k_1,k_2}$ is precisely
We define a set function 
\begin{equation*}
\nu_y\left(E_{k_1,k_2}\right) = E(x;a)^{-1}|E_{k_1,k_2}|.
%\left|\left\{ n \leq x : \mu^2(n(n+1)) = 1, k_1|n, k_2|(n+d), (n,P/k_1) = 1, (n+1,P/k_2)=1\right\}\right|.
\end{equation*}
Note that the sets $\{E_{k_1,k_2}\}_{k_1,k_2|P \atop (k_1,k_2) = 1}$ are mutually disjoint, and their union is precisely the set of all $n \leq x$ coprime to $a$ such that $n,n+a \in S$.
%and, by the sieve,
%\begin{align*}
%\left|E(x)-\left|\bigcup_{2 \leq k_1,k_2|P \atop (k_1,k_2) =1} E_{k_1,k_2}\right|\right| &\leq \sum_{j = 0,1} \left|\left\{n \leq x : \mu^2(n(n+1)) = 1, k \nmid (n+j) \forall k|P\right\}\right| \\
%&\ll \left|\left\{ n \leq x : P^-(n) > z\right\}\right| \ll \frac{x}{\log z}.
%\end{align*}
%Thus, $\nu\left(\bigcup_{2 \leq k_1,k_2|P \atop (k_1,k_2) = 1}\right) = 1 + O\left(\frac{1}{\log z}\right)$. Our approximation to $\nu$ will be a probability measure. \\
We use a sieve to estimate each individual $\nu_y(E_{k_1,k_2})$. 
\begin{lem}\label{ROSSER}
Let $x \geq 3$ and $y \leq D_1,D_2 < x$. Put $s_j := \log D_j/\log y$, for $j = 1,2$. Then if $k_1,k_2 | P$ are coprime, 
%and $D_1D_2k_1k_2 < x^{\frac{1}{2}}$ then
\begin{align*}
\nu_y\left(E_{k_1,k_2}\right) &= f(k_1)f(k_2)\prod_{p | P/k_1k_2} (1-2f(p)) \left(1+O\left(\max_{j = 1,2}\left\{e^{-(1+o(1))s_j \log s_j} + \log^{-\frac{1}{6}} D_j\right\}\right)\right) \\
&+O\left(\sum_{m_1m_2 \leq k_1k_2D_1D_2} f(k_1)f(k_2)|R(x; m_1,m_2)|\right).
%&+ O_{\e}\left(\left(\frac{x}{k_1k_2}\right)^{\frac{2}{3}+\e} (D_1D_2)^{\frac{1}{3}-\e}\right).
\end{align*}
\end{lem}
\begin{proof}
For $j = 1,2$, let $\Lambda_j^{+} := \{\lambda_j^+(d)\}_{d|P, d \leq D_j}$ and $\Lambda_j^- := \{\lambda_j^-(d)\}_{d | P , d\leq D_j}$ be upper and lower bound sieves, i.e., sequences of real numbers that satisfy 
\begin{equation*}
1 \ast \lambda_j^- (d) \leq 1 \ast \mu(d) \leq 1 \ast \lambda_j^+(d)
\end{equation*}
for each $d | P$ with $d \leq D_j$. We take $\lambda_j^+$ and $\lambda_j^-$ to be the upper and lower bound sieve weights of Rosser-Iwaniec, respectively (see Chapter 11 of \cite{FI}). \\
We can compose these sieve weights to produce upper and lower bound sieves by defining the two-dimensional weights
\begin{align*}
\varphi^+(m,n) &:= \left(1\ast \lambda_1^+\right)(m)\left(1\ast \lambda_2^+\right)(n)\\
\varphi^-(m,n) &:= \left(1\ast \lambda_1^-\right)(m)\left(1\ast \lambda_2^+\right)(n) + \left(1\ast \lambda_1^+\right)(m)\left(1\ast \lambda_2^-\right)(n) - \left(1\ast \lambda_1^+\right)(m)\left(1\ast \lambda_2^+\right)(n).
\end{align*}
That $\varphi^+$ is an upper bound sieve weight is immediate; that $\varphi^-$ is a lower bound sieve weight follows from the identity 
\begin{equation*}
\varphi^-(m,n) = \left(1\ast \lambda_1^-\right)(m)\left(1\ast \lambda_2^-\right)(n) - \left(\left(1\ast \lambda_1^+\right)(m)-\left(1\ast \lambda_1^-\right)(m)\right)\left(\left(1\ast \lambda_2^+\right)(n)-\left(1\ast \lambda_2^-\right)(n)\right).
\end{equation*}
Therefore, by Definition \ref{SIFTABLE} with $X = E(x;a)$,
\begin{align}
|E_{k_1,k_2}| &\leq X\mathop{\sum_{d_1|P/k_1}\sum_{d_2|P/k_2}}_{(k_1d_1,k_2d_2) = 1} \lambda_1^+(d_1)\lambda_2^+(d_2)f(k_1k_2d_1d_2) + \mathop{\sum_{d_1|P/k_1}\sum_{d_2|P/k_2}}_{(k_1d_1,k_2d_2) = 1} f(k_1k_2)|R(x;k_1d_1,k_2d_2)| \label{UPPER} \\
|E_{k_1,k_2}| &\geq X\mathop{\sum_{d_1|P/k_1}\sum_{d_2|P/k_2}}_{(k_1d_1,k_2d_2) = 1}\left(\lambda_1^+(d_1)\lambda_2^-(d_2) + \lambda_1^-(d_1)\lambda_2^+(d_2) - \lambda_1^+(d_1)\lambda_2^+(d_2)\right)f(k_1k_2d_1d_2) \\
&- \mathop{\sum_{d_1|P/k_1}\sum_{d_2|P/k_2}}_{(k_1d_1,k_2d_2) = 1} f(k_1k_2)|R(x;k_1d_1,k_2d_2)| \label{LOWER}.
\end{align}
Since $f$ is multiplicative, for any of the sign pairs $(\eta_1,\eta_2) \in \{(+,+),(+,-),(-,+)\}$,
\begin{equation*}
\mathop{\sum_{d_1|P/k_1}\sum_{d_2|P/k_2}}_{(k_1d_1,k_2d_2) = 1} \lambda_1^{\eta_1}(d_1)\lambda_2^{\eta_2}(d_2)f(k_1k_2d_1d_2) = f(k_1)f(k_2)\sum_{d_1|P/k_1 \atop (d_1,k_2) = 1}\lambda_1^{\eta_1}(d_1)f(d_1)\sum_{d_2|P/k_2 \atop (d_2,k_1d_1) = 1} \lambda_2^{\eta_2}(d_2)f(d_2).
\end{equation*}
Put $h_{k}(m) := f(m)1_{(m,k) = 1}$. Then by Theorem 11.12 of \cite{FI}, we have
\begin{align*}
&\sum_{d_1|P/k_1 \atop (d_1,k_2) = 1}\lambda_1^{\eta_1}(d_1)f(d_1)\sum_{d_2|P/k_2 \atop (d_2,k_1d_1) = 1} \lambda_2^{\eta_2}(d_2)f(d_2) = \sum_{d_1|P/k_1 }\lambda_1^{\eta_1}(d_1)h_{k_2}(d_1)\sum_{d_2|P} \lambda_2^{\eta_2}(d_2)h_{k_1k_2d_1}(d_2) \\
&= \left(1+O\left(e^{-s_2\log s_2} + \log^{-\frac{1}{6}} D_2\right)\right) \sum_{d_1|P/k_1}\lambda_1^{\eta_1}(d_1)h_{k_2}(d_1)\prod_{p|P/(k_1k_2d_1)}\left(1-f(p)\right) \\
&= \left(1+O\left(e^{-(1+o(1))s_2\log s_2} + \log^{-\frac{1}{6}} D_2\right)\right) \prod_{p|P/k_1k_2}(1-f(p))\sum_{d_1|P/k_1 }\lambda_1^{\eta_1}(d_1)h_{k_2,k_1}(d_1)\prod_{p|d_1}\left(1-f(p)\right)^{-1} \\
&= \left(1+O\left(\max_{j = 1,2} \left\{e^{-(1+o(1))s_j\log s_j} + \log^{-\frac{1}{6}} D_j\right\}\right)\right) \prod_{p|P/k_1k_2}(1-f(p))\left(1-f(p)(1-f(p))^{-1}\right) \\
&= \left(1+O\left(\max_{j = 1,2} \left\{e^{-(1+o(1))s_j\log s_j} + \log^{-\frac{1}{6}} D_j\right\}\right)\right) \prod_{p|P/k_1k_2}(1-2f(p)).
\end{align*}
This last estimate together with \eqref{UPPER} and \eqref{LOWER} gives
\begin{align*}
|E_{k_1,k_2}| &= \left(1+O\left(\max_{j = 1,2} \left\{e^{-(1+o(1))s_j\log s_j} + \log^{-\frac{1}{6}} D_j\right\}\right)\right) Xf(k_1)f(k_2)\prod_{p|P/k_1k_2}(1-2f(p)) \\
&+ O\left(\sum_{d_1 \leq D_1, d_2 \leq D_2} f(k_1k_2)|R(x;k_1d_1,k_2d_2)|\right).
\end{align*}
%and it remains to estimate the contribution from the remainder term.  Since $S$ is siftable,
Moreover, we can trivially bound the remainder by
\begin{equation*}
\sum_{d_1 \leq D_1, d_2 \leq D_2} |R(x;k_1d_1,k_2d_2)| \leq \sum_{m_1m_2 \leq D_1D_2k_1k_2} |R(x;m_1,m_2)|.
% \ll \frac{x}{(\log x)^3}.
\end{equation*}
%By Lemma \ref{SIEVESTUFF}, we have $R_{k_1,k_2}(d_1,d_2) \ll_{\e} (x/(k_1k_2d_1d_2))^{\frac{2}{3}+\e}$ for $k_1k_2 < x$ and $d_j \leq D_j$, so that
%\begin{equation*}
%\sum_{d_1 \leq D_1, d_2 \leq D_2} |R(x;k_1d_1,k_2d_2)| \ll_{\e} \left(\frac{x}{k_1k_2}\right)^{\frac{2}{3}+\e} \prod_{j = 1,2} \left(\sum_{d_j \leq D_j} d_j^{-\left(\frac{2}{3}+\e\right)}\right) \ll_{\e} \left(\frac{x}{k_1k_2}\right)^{\frac{2}{3}+\e} (D_1D_2)^{\frac{1}{3}-\e}.
%\end{equation*}
This completes the proof.
\end{proof}
%Specifying our parameters in a manner that will be suitable in the sequel, gives the following formulation.
\begin{cor} \label{ROSSCOR}
Fix $k_1,k_2|P$ with $k_1k_2 \leq x^{\frac{\theta}{2}}$, where $\theta$ is associated to $S$. Let $y = x^{\frac{1}{\beta}}$ and $D := \left(\frac{x}{k_1k_2}\right)^{\frac{\theta}{4}}$. Then
\begin{equation*}
\nu_y(E_{k_1,k_2}) = \left(1+O\left(e^{-\frac{\theta}{8}\beta \log \beta} + \log^{-\frac{1}{6}} x\right)\right)f(k_1)f(k_2)\prod_{p|P/k_1k_2}(1-2f(p)).
\end{equation*}
\end{cor}
\begin{proof}
Choose $D_1 = D_2 = D$ in Lemma \ref{ROSSER} so that $s_1 = s_2 \geq \theta \beta/8$. Note that $k_1k_2D_1D_2 = x^{\frac{\theta}{2}}(k_1k_2)^{1-\frac{\theta}{4}} \leq x^{\theta}$. Thus, since $S$ is siftable, the sum of the remainder terms in Lemma \ref{ROSSER} can be estimated as
\begin{equation*}
\sum_{m_1m_2 \leq k_1k_2D_1D_2} |R(x;m_1,m_2)| \leq \sum_{m_1m_2 \leq x^{\theta}} |R(x;m_1,m_2)| \ll \frac{x}{\log^3 x}.
\end{equation*}
Note moreover that 
\begin{equation*}
\prod_{p |P/k_1k_2}\left(1-2f(p)\right) \gg \prod_{p \leq y} \left(1-\frac{1}{p}\right)^2 \gg (\log x)^{-2},
\end{equation*}
so that $x(\log x)^{-3} \ll (\log x)^{-1}\prod_{p | P/k_1k_2}\left(1-2f(p)\right)$. The lemma now follows from Lemma \ref{ROSSER} upon dividing by $E(x;a)$, since $E(x;a) \gg x$ by assumption.
\end{proof}
Let $\mathfrak{S}$ denote the algebra of subsets of $\mb{N} \cap [1,x]$ generated by the sets $E_{k_1,k_2}$, with $k_1$ and $k_2$ coprime divisors of $P$. As mentioned, the sets $E_{k_1,k_2}$ are mutually disjoint. Thus, any $A \in \mathfrak{S}$ can be written as a $$A = \bigcup_{(k_1,k_2) \in C} E_{k_1,k_2},$$ where $C$ is a set of pairs of coprime divisors of $P$. Define
\begin{align*}
\nu_y(A) &= \sum_{(k_1,k_2) \in C} \nu_y(E_{k_1,k_2}) \\
&=E(x;a)^{-1}\sum_{(k_1,k_2) \in C} \sum_{n \leq x \atop (n,P) = k_1, (n+a,P) = k_2} 1_S(n)1_S(n+a)1_{(n,a) = 1}.
\end{align*}
We construct a measure $\sigma_y$ whose value, determined on the sets $E_{k_1,k_2}$ generating $\mathfrak{S}$, is defined by $$\sigma_y(E_{k_1,k_2}) = f(k_1)f(k_2)\prod_{p|P/k_1k_2}(1-2f(p)),$$ where $f$ is the multiplicative function in \eqref{SIEVE}. By Corollary \ref{ROSSCOR}, whenever $k_1k_2 \leq x^{\frac{\theta}{2}}$ we have
\begin{equation}\label{COMP}
\nu_y(E_{k_1,k_2}) = \left(1+O\left(e^{-\frac{1}{6}\beta \log \beta} + \log^{-\frac{1}{6}} x\right)\right) \sigma_y(E_{k_1,k_2}).
% + O\left(x^{-\frac{1}{24}}\right).
\end{equation}
%In the one-dimensional case, this procedure is due originally to Kubilius. 
The following lemma is inspired by the one-dimensional treatment of Kubilius' model given in Chapter 3 of \cite{Ell}. \\
Put $R(x,y) := e^{-\frac{\theta}{8}\min\{\theta,\frac{1}{2}\}\beta \log \beta} + \log^{-\frac{1}{6}} x$. 
\begin{lem} \label{TVONE}
As defined, $\sigma_y$ is a probability measure. Moreover, for any $A \in \mc{S}$ we have 
\begin{equation*}
\nu_y(A) = \left(1+O\left(R(x,y)\right)\right)\sigma_y(A) + O\left(R(x,y)\right).
%+ O\left(\frac{1}{\log_2 x}\right).
\end{equation*}
\end{lem}
\begin{proof}
Recall that the collection $\{E_{k_1,k_2}\}_{k_1,k_2|P \atop (k_1,k_2) = 1}$ partitions the set of $n \leq x$ such that $n,n+a \in S$ and $(n,a) = 1$. Thus, to establish that $\sigma_y$ is a probability measure it suffices to show that
\begin{equation} \label{PROB}
\sum_{k_1,k_2|P \atop (k_1,k_2) = 1} \sg_y(E_{k_1,k_2}) = \sum_{k_1,k_2|P \atop (k_1,k_2) = 1} f(k_1)f(k_2)\prod_{p|P/k_1k_2}(1-2f(p)) = 1.
\end{equation}
This is a straightforward computation. Indeed, as $P$ is squarefree,
\begin{align*}
&\sum_{k_1,k_2|P \atop (k_1,k_2) = 1} f(k_1)f(k_2)\prod_{p|P/k_1k_2}(1-2f(p)) \\
&= \prod_{p|P}\left(1-2f(p)\right)\sum_{k_1|P}f(k_1)\prod_{p|k_1}(1-2f(p))^{-1}\sum_{k_2|P/k_1}f(k_2)\prod_{p|k_2}(1-2f(p))^{-1} \\
&= \prod_{p|P}\left(1-2f(p)\right)\sum_{k_1|P}f(k_1)\prod_{p|k_1}(1-2f(p))^{-1}\prod_{p|P/k_1}\left(1+f(p)(1-2f(p))^{-1}\right) \\
&= \prod_{p|P}\left(1-2f(p)\right)\left(1+f(p)(1-2f(p))^{-1}\right)\sum_{k_1|P}f(k_1)\prod_{p|k_1}(1-2f(p))^{-1}\left(1+f(p)(1-2f(p))^{-1}\right)^{-1} \\
&= \prod_{p|P}\left(1-f(p)\right)\left(1+f(p)(1-f(p))^{-1}\right) = 1.
\end{align*}
We next establish that
\begin{equation} \label{TAIL}
\sum_{k_1,k_2|P \atop (k_1,k_2) = 1, k_1k_2 > x^{\frac{\theta}{2}}} \sigma_y(E_{k_1,k_2}) \ll R(x,y).
\end{equation}
Assume that \eqref{TAIL} has been proven. Note then that since $\nu_y$ is a probability measure, \eqref{COMP} and \eqref{PROB} together imply that
\begin{align*}
\sum_{k_1,k_2|P \atop (k_1,k_2) = 1, k_1k_2 > x^{\frac{\theta}{2}}} \nu_y(E_{k_1,k_2}) &= 1-\sum_{k_1,k_2|P \atop (k_1,k_2) = 1, k_1k_2 \leq x^{\frac{\theta}{2}}} \nu_y(E_{k_1,k_2}) \\
&= 1-\left(1+O\left(R(x,y)\right)\right) \sum_{k_1,k_2|P \atop (k_1,k_2) = 1, k_1k_2 \leq x^{\frac{\theta}{2}}} \sigma_y(E_{k_1,k_2}) \\
&= \sum_{k_1,k_2|P \atop (k_1,k_2) = 1, k_1k_2> x^{\frac{\theta}{2}}} \sigma_y(E_{k_1,k_2}) + O\left(R(x,y)\right) = O\left(R(x,y)\right).
\end{align*}
Thus, if $A = \bigcup_{(k_1,k_2) \in C} E_{k_1,k_2}$ and we write $$C^+ := C \cap \{k_1,k_2| P : (k_1,k_2) = 1, k_1k_2 > x^{\frac{\theta}{2}}\}$$ then by \eqref{TAIL},
\begin{align*}
\nu_y(A) &= \sum_{(k_1,k_2) \in C} \nu_y(E_{k_1,k_2}) \\
&= \left(1+O\left(R(x,y)\right)\right)\sum_{(k_1,k_2) \in C} \sigma_y(E_{k_1,k_2}) + \sum_{(k_1,k_2) \in C^+}\left(\nu(E_{k_1,k_2})-\left(1+O\left(R(x,y)\right)\right)\sigma_y(E_{k_1,k_2})\right) \\
&= \left(1+O\left(R(x,y)\right)\right)\sigma_y(A) + O\left(\sum_{k_1,k_2|P \atop (k_1,k_2) = 1, k_1k_2 > x^{\frac{\theta}{2}}}(\sigma_y(E_{k_1,k_2})+\nu_y(E_{k_1,k_2}))\right) \\
&= \left(1+O\left(R(x,y)\right)\right)\sigma_y(A) + O\left(R(x,y)\right),
\end{align*}
as claimed. Thus, it remains to prove \eqref{TAIL}. \\
To this end we use Rankin's trick. Let $\e > 0$ be sufficiently small but fixed. First, by assumption $f(p) \leq \frac{1}{p}$ when $p \geq 2$. Thus, 
\begin{equation*}
f(p)^{1-\e} - f(p) \leq \frac{1}{p}\left(p^{\e}-1\right) \ll \e y^{\e} \frac{\log p}{p},
\end{equation*}
for each $p \leq y$ once $\e < 1/2$.  Hence,
%Furthermore, observe that 
%\begin{equation*}
%(1-2f(p))^{-1} \leq (1-2/p)^{-1} \leq 1+\frac{3}{p}
%\end{equation*}
%for $p \geq 7$. Thus,
\begin{align*}
\sum_{k_1,k_2|P \atop (k_1,k_2) = 1, k_1k_2 > x^{\frac{\theta}{2}}} \sigma_y(E_{k_1,k_2}) &= \prod_{p \leq y \atop p \nmid a} \left(1-2f(p)\right) \sum_{k_1,k_2|P \atop (k_1,k_2) = 1, k_1k_2 > x^{\frac{\theta}{2}}} f(k_1)f(k_2) \prod_{q |k_1k_2} \left(1-2f(q)\right) \\
&\leq \prod_{p \leq y \atop p \nmid a} \left(1-2f(p)\right)x^{-\frac{\e\theta}{2}} \left(\sum_{k|P} f(k)^{1-\e}\prod_{q|k}(1-2f(q))^{-1}\right)^2 \\
&\ll_{\e} \prod_{p \leq y \atop p \nmid a} \left(1-2f(p)\right)x^{-\frac{\e\theta}{2}} \exp\left(2\sum_{q \leq y \atop q \nmid a} f(q)^{1-\e}\left(1-2f(p)\right)^{-1}\right) \\
&\ll x^{-\frac{\e\theta}{2}} \exp\left(2\sum_{p \leq y} \left(f(p)^{1-\e} - f(p)\right)\right)\\
&\leq \exp\left(\e\left(2y^{\e}\sum_{p \leq y} \frac{\log p}{p}-\frac{\theta}{2}\log x\right)\right).
\end{align*}
%Now note that $f(p) \geq \frac{1}{p}(1-1/p)$, so that
%\begin{equation*}
%-\sum_{p \leq z} \log\left(1-2f(p)\right) \geq O(1) + 2\sum_{p \leq z} \frac{1}{p}\left(1-\frac{1}{p}\right) = 2\log_2 z + O(1).
%\end{equation*}
By Mertens' theorem, the exponential is $\ll \left(y^{\e}\right)^{2y^{\e}}x^{-\frac{\e\theta}{2}}$. 
Choosing $\e := \frac{C}{\log y}$ for $C > 0$ to be chosen, it follows that
\begin{equation*}
\sum_{k_1,k_2|P \atop (k_1,k_2) = 1, k_1k_2 > x^{\frac{\theta}{2}}} \sg_y(E_{k_1,k_2}) \ll \exp\left(2C\left(e^{C}-\frac{\theta\log x}{4 \log y}\right)\right) = \exp\left(2C\left(e^C - \frac{\theta}{4}\beta\right)\right).
\end{equation*}
Put $\lambda := \frac{\theta}{4}\beta$. Choosing $C := \log\lambda - \log_2 \lambda$ yields
\begin{equation*}
\sum_{k_1,k_2|P \atop (k_1,k_2) = 1, k_1k_2 > x^{\frac{\theta}{2}}} \sigma_y(E_{k_1,k_2})  \ll e^{-\frac{2}{3}\lambda \log \lambda} \ll e^{-\frac{\theta}{6}\beta \log \beta}
\end{equation*}
and \eqref{TAIL} follows.
%Using our choice of $y$ finally proves \eqref{TAIL} for $x$ sufficiently large.
\end{proof}
For each $p \leq y$ let $X_p := (\tilde{\theta}_p,\tilde{\theta}_p'): [1,x]^2 \ra \{0,1\}$ denote the independent Bernoulli random vector satisfying
\begin{equation*}
X_p = \begin{cases} (1,0) &\text{ Prob} = f(p) \\ (0,1) &\text{ Prob} = f(p) \\ (0,0) &\text{ Prob} = 1-2f(p).\end{cases}
\end{equation*}
Put $\Sigma(y) := \sum_{p \leq y} X_p$, and note that $\sigma_y$ is the law of $\Sigma(y)$. Indeed, we can identify the event $\{\Sigma(y) = (l_1,l_2)\}$ with the set
\begin{equation*}
\mc{S}(l_1,l_2) := \bigcup_{k_1,k_2 | P \atop \omega(k_j) = l_j}^{\ast} E_{k_1,k_2},
\end{equation*}
where the asterisk indicates a union over $(k_1,k_2)= 1$. Since $X_p$ and $X_q$ are independent in the probability space $([1,x], \mathfrak{S},\sigma_y)$, 
\begin{align*}
\sigma_y(E_{k_1,k_2}) &= \prod_{p|k_1k_2} f(p) \prod_{q |P/k_1k_2} (1-2f(q)) \\
&= \prod_{p|k_1} \mb{P}(X_p = (1,0)) \prod_{q|k_2} \mb{P}(X_q = (0,1)) \prod_{r|P/k_1k_2} \mb{P}(X_r = (0,0)) \\
&= \mb{P}\left(X_p = (1,0) \text{ iff } p|k_1, X_q = (0,1) \text{ iff } q|k_2\right) = \mb{P}\left(\tilde{\theta}_p = 1 \text{ iff } p|k_1, \tilde{\theta}_p' = 1 \text{ iff } p|k_2\right).
\end{align*}
Therefore, it follows that
\begin{align*}
\sigma_y(\mc{S}(l_1,l_2)) &= \sum_{k_1,k_2|P \atop \omega(k_j) = l_j}^{\ast} \sigma_y(E_{k_1,k_2}) = \sum_{k_1,k_2|P \atop \omega(k_j) = l_j}^{\ast} \mb{P}\left(\tilde{\theta}_p = 1 \text{ iff } p|k_1, \tilde{\theta}_p' = 1 \text{ iff } p|k_2\right) = \mb{P}(\Sigma(y) = (l_1,l_2)).
\end{align*}
We note that since $\sigma_y$ is an atomic measure, the relationship between $\sigma_y$ and $\nu_y$ in Lemma \ref{TVONE} continues to hold when both $\sum_{p \leq y} X_p$, and $(\omega_y,\omega_y \circ \sg^a)$ are respectively centred and normalized by $\sum_{p \leq y} f(p)$. \\
%and $(\omega,\omega \circ S)$ is transformed into $(\tilde{\omega}_x,\tilde{\omega}_x \circ S)$. \\
%, and suppose that for $p,q \leq z$ distinct, $\theta_p$ and $\theta_q$ are independent. Let $\{\theta_p'\}_{p \leq z}$ be an i.i.d. copy of the sequence $\{\theta_p\}_{p\leq z}$.  We claim that $\mu$ is the law of the random vector $(\hat{\omega},\hat{\omega}')$. 
%We will also show that the total variation norm of the law $\nu$ of the deterministic vector $(\omega,\omega \circ S)$ relative to $\mu$ is small. The relevance of this intermediate step stems from the following result, due to Roos. \\
To proceed, we will require the following result due to Roos \cite{ROOS}. In what follows, fix $d \in \mb{N}$ and for each $1 \leq j \leq d$ let $\mbf{e}_j \in \mb{R}^d$ be the unit vector whose only non-zero component is in the $j$th position.
\begin{lem}[\cite{ROOS}, Corollary 1]\label{ROOSTHM}
Let $\{X_k\}_{k \geq 1}$ be a sequence of Bernoulli random $d$-vectors. 
%and write $\mbf{X}_k = (X_{k,j})_{1 \leq j \leq d$. 
Let $\mb{P}(X_k = \mbf{e}_j) = p_{k,j}$, $p_k := \sum_{1 \leq j \leq k} p_{k,j}$, and suppose that $\mb{P}(X_k = \mbf{0}) = 1-p_k$. Put $S_n := \sum_{1 \leq k \leq n} X_k$ and $\lambda_n(j):= \sum_{1 \leq k \leq n} p_{k,j}$. Then if $Z \sim \text{Poi}(\lambda_n(1),\ldots\lambda_n(d))$ with pairwise independent coordinates,
\begin{equation*}
d_{TV}(\mc{L}(S_n),\mc{L}(Z)) \ll \sum_{j \leq d}\min\left\{1,\lambda_n(j)^{-1}\right\} \sum_{1 \leq k \leq n} p_{k,j}^2.
\end{equation*}
\end{lem}
An immediate consequence of our foregoing analysis and Lemma \ref{ROOSTHM} is the following.
\begin{lem} \label{TVPOISSON}
Let $Z = (Z_1,Z_2) \sim \text{Poi}(\lambda(y),\lambda(y))$, where $\lambda(y) := \sum_{p \leq y} f(p)$. Then
\begin{equation*}
d_{TV}(\nu_y,\mc{L}(Z)) \ll \frac{1}{\lambda(y)} + R(x,y).
\end{equation*}
\end{lem}
\begin{proof}
By the triangle inequality, we have $d_{TV}(\nu_y,\mc{L}(Z)) \leq d_{TV}(\nu_y,\sigma_y) + d_{TV}(\sigma_y,\mc{L}(Z))$. Lemma \ref{TVONE} implies that $d_{TV}(\nu_y,\sigma_y) \ll R(x,y)$. We apply Lemma \ref{ROOSTHM} with $n = \pi(y)$ and $p_{q,1} = p_{q,2} = f(q)$ for each $q\leq y$. Thus, $\lambda_n(1) = \lambda_n(2) = \sum_{q \leq y} f(q)$. Moreover, as $f(q) \leq 1/q$ for each prime $q$, $\sum_{q \leq y} p_{q,j}^2 \ll 1$ for $j = 1,2$. Thus, Lemma \ref{ROOSTHM} indeed implies that $d_{TV}(\sigma_y,\mc{L}(Z)) \ll \lambda(y)^{-1}$, and the proof is complete. %Finally, as $\log_2 z = \log_2 x - \log_4 x + O(1)$, we have $\log_2^{-1} z \ll \log_2^{-1} x$, and the lemma is proved.
\end{proof}
\section{Uniform Approximation of $\phi_x$ by Normal and Poisson Characteristic Functions}
In this section, we show that $\phi_{x,y}$ uniformly approximates the cahracteristic function of $\Phi_{(2)}$. We use this data to establish and estimate for the $L^{\infty}$ distance between $F_{x,S,a}$ and $\Phi_{(2)}$. The tool that makes this connection possible is the following two-dimensional version of a smoothing lemma due to Ess\'{e}en (see Chapter XVI, Section 6 in \cite{FEL}).
\begin{lem}\label{ESSEEN}
Let $T \geq 1$. Let $G: \mb{R}^2 \ra \mb{R}$ be a differentiable function with $\|\nabla G\|_{\infty} \leq 1$, and such that $G(t,t') \ra 1$ as $t,t' \ra \infty$ and $G(t,t') \ra 0$ as $t,t' \ra -\infty$. Let $\chi$ be the Fourier transform of the Lebesgue-Stieltjes measure $dG$. Let $H$ be the distribution function of a bivariate random vector whose characteristic function is $\phi$. Furthermore, suppose that for any fixed $\mbf{z}$,
\begin{equation} \label{TAILS}
\left|G(\mbf{z}) - H(\mbf{z}) - \int_{-T}^{z_1}\int_{-T}^{z_2} d(G-H)(\mbf{w})\right| \ll T^{-1}.
\end{equation}
Then 
%Assume, moreover, that $G-H \in L^1(\mb{R}^2)$.  Then
\begin{equation} \label{SMOOTHESS}
\|G-H\|_{L^{\infty}(\mb{R}^2)} \ll \int_{\mc{R}_T} \frac{\left|\phi(t,t')-\chi(t,t')\right|}{|tt'|}dtdt' + T^{-1},
\end{equation}
where $\mc{R}_T := \{\mbf{u} \in [-T,T]^2 : T^{-3} < |u_1|,|u_2| \leq T\}$.
\end{lem}
The proof proceeds along similar lines to Ess\'{e}en's one-dimensional lemma, but the author did not find it proved in the context of bivariate distribution functions. Thus, we prove it here for completeness.
\begin{proof}
Put $\Delta := H-G$ and let $\mbf{z}_0 \in \mb{R}^2$ be such that $\Delta(\mbf{z}_0) = \|\Delta\|_{L^{\infty}(\mb{R}^2)}$ (we may assume positivity by defining $\Delta$ as $G-H$ instead, if necessary). Put $w_T(y) := \frac{1-\cos(Ty)}{\pi y^2}$ and let $W_T(u,v) := w_T(u)w_T(v)$. Note that 
\begin{equation*}
\hat{W}(r,s) = \hat{w}(r)\hat{w}(s) = \begin{cases} T^2\left(1-\frac{|r|}{T}\right)\left(1-\frac{|s|}{T}\right) & \text{ if $\max\{|r|,|s|\} \leq T$} \\ 0& \text{ otherwise}. \end{cases}
\end{equation*}
Also, write $\Delta_T := \Delta \ast W_T$. Let $h_1,h_2 \geq 0$ and $\mbf{h} := (h_1,h_2)$. Note that since $H$ is a distribution function, $H(\mbf{z}_0) \geq H(\mbf{z}_0- \mbf{h})$. Thus, by Taylor's theorem,
\begin{equation*}
\Delta(\mbf{z_0}) - \Delta(\mbf{z_0}-\mbf{h}) \geq H(\mbf{z_0}) - H(\mbf{z_0}-\mbf{h}) -|G(\mbf{z_0}-\mbf{h})-G(\mbf{z_0})| \geq -|(\nabla G)(\mbf{h}) \cdot \mbf{h}| \geq -(h_1+h_2)
\end{equation*}
Now, let $\delta := \Delta(\mbf{z_0})$, let $B := B_{\delta}(\mbf{z_0})$ be the $\delta$-ball (with respect to the $L^{\infty}$-norm), and put $\mbf{z} := \mbf{z_0} - \delta\cdot (1,1)$. Then for all $\mbf{u} \in B$, we have $H(\mbf{z}+\mbf{u}) \leq H(\mbf{z_0})$. Noting that $w_T$ is even, we have
\begin{align*}
\int_{-\delta}^{\delta} \int_{-\delta}^{\delta} \Delta(\mbf{z}+\mbf{h}) W_T(\mbf{h}) d\mbf{h} &\leq \int_{-\delta}^{\delta}\int_{-\delta}^{\delta} \left(\Delta(\mbf{z_0}) -2\delta + h_1+h_2\right)W_T(\mbf{h}) d\mbf{h} \\
&= -4\delta \left(\int_0^{\delta} w_T(u) du\right)^2 = -4\delta T^2\left(\int_0^{T\delta} \frac{1-\cos u}{\pi u^2} du\right);
\end{align*}
this last statement follows by definition of $\delta$ and because $h_jw_T(h_j)$ is odd. Moreover,
\begin{align*}
\int_{\|\mbf{h}\|_{\infty} > \delta} \Delta(\mbf{z}+\mbf{h})W_T(\mbf{h}) &\leq 4\delta \left(\int_0^{\infty} \frac{1-\cos(Tx)}{\pi x^2} dx\right)\int_{|t| > \delta} \frac{1-\cos(Tx)}{\pi x^2} \\
&= 4\delta T^2\int_{u > T\delta} \frac{1-\cos u}{\pi u^2} du \leq \frac{4}{\pi^2}T.
\end{align*}
Thus, we have
\begin{align}
|\Delta_T(\mbf{z})| &\geq \left|\int_{-\delta}^{\delta}\int_{-\delta}^{\delta} \Delta(\mbf{z}+\mbf{h})W_T(\mbf{h}) d\mbf{h}\right| -\frac{4}{\pi^2}T \nonumber\\
&\geq 4T^2 \left(\delta \left(\int_0^{T\delta} \frac{1-\cos u}{\pi u^2} du \right)^2 - \frac{1}{\pi^2 T}\right) \\
&\geq 4T^2 \left(\delta - \frac{3}{\pi^2 T}\right) = 4T^2\left(\|\Delta\|_{\infty}-\frac{3}{\pi^2 T}\right). \label{LEFT}
\end{align}
It remains to show that $|\Delta_T(\mbf{z})|$ is bounded by the integral on the right side of \eqref{SMOOTHESS}. \\
Now, let $d\mu_T := W_T d\lambda$, where $d\lambda$ is Lebesgue measure on $\mb{R}^2$. Let $B(\mbf{z})$ denote the infinite box $(-\infty,z_1] \times (-\infty,z_2]$. Then the convolution of measures $\Delta \ast \mu_T$ can be written as
\begin{equation*}
\Delta \ast \mu_T(B(\mbf{z})) 
%= \int_{\mb{R}^2} \int_{\mb{R}^2} 1_{B(\mbf{z})}(\mbf{u} + \mbf{v}) d\Delta(\mbf{u}) d\mu_T(\mbf{v}) 
= \int_{\mb{R}^2} \Delta(\mbf{z} + \mbf{u}) d\mu_T(\mbf{v}) = \int_{\mb{R}^2} \Delta(\mbf{z} + \mbf{u}) W_T(\mbf{v}) d\mbf{v} = \Delta_T.
\end{equation*}
Now, by Fourier inversion, we have
\begin{align*}
\Delta\ast \mu_T(B(\mbf{z})) &= \int_{-\infty}^{z_1}\int_{-\infty}^{z_2} d(\Delta \ast \mu_T)(\mbf{w}) \\
&= \frac{1}{2\pi} \int_{-\infty}^{z_1}\int_{-\infty}^{z_2} d\xi_1d\xi_2\left(\int_{\mb{R}^2} e^{-i\mbf{u} \cdot \mbf{\xi}} \hat{w}_T(u_1) \hat{w}_T(u_2)(\phi(\mbf{u})-\chi(\mbf{u})) du_1 du_2\right) \\
&= \frac{T^2}{2\pi} \lim_{\tau_1,\tau_2 \ra \infty} \left(\int_{-\tau_1}^{z_1}\int_{-\tau_2}^{z_2} d\mbf{\xi}\left(\int_{[-T,T]^2} e^{-i\mbf{u} \cdot \mbf{\xi}} \left(1-\frac{|u_1|}{T}\right)\left(1-\frac{|u_2|}{T}\right)(\phi(\mbf{u})-\chi(\mbf{u})) d\mbf{u}\right)\right). 
%= \frac{1}{2\pi} \int_{[-T,T]^2} e^{-i(uz_1 + vz_2)} \left(1-\frac{|u|}{T}\right)\left(1-\frac{|v|}{T}\right)(\phi(u,v)-\chi(u,v)) du dv
\end{align*}
This shows that the Radon-Nikodym derivative of $\Delta \ast \mu_T$ with respect to Lebesgue measure is precisely
\begin{equation*}
h_{T}(\mbf{\xi}) := \frac{T^2}{2\pi} \int_{[-T,T]^2} e^{-i\mbf{u}\cdot \mbf{\xi}} \left(1-\frac{|u_1|}{T}\right)\left(1-\frac{|u_2|}{T}\right)(\phi(\mbf{u})-\chi(\mbf{u})) du_1 du_2.
\end{equation*}
For $\e > 0$, let $\mc{S}_{\e} := \{\mbf{u} \in [-T,T]^2 : \e \leq |u_1|,|u_2| \leq T\}$ and choose a smooth function $k_{\e}: \mb{R}^2 \ra [0,1]$ with support lying inside of $\mc{S}_{\e}$ such that $k_{\e} \ra 1$ uniformly on $[-T,T]^2$ as $\e \ra 0$. Define
\begin{equation*}
h_{T,\e}(\mbf{\xi}) := \frac{T^2}{2\pi} \int_{[-T,T]^2} e^{-i\mbf{u}\cdot \mbf{\xi}} \left(1-\frac{|u_1|}{T}\right)\left(1-\frac{|u_2|}{T}\right)(\phi(\mbf{u})-\chi(\mbf{u}))k_{\e}(\mbf{u}) du_1 du_2.
\end{equation*}
By the Dominated Convergence theorem, $h_{T,\e} \ra h_T$ as $\e \ra 0$ uniformly on $[-T,T]^2$. Now, define the functions
\begin{align*}
J_{\e,1}(\mbf{\xi}) &:= \frac{T^2}{2\pi} \int_{\mc{S}_{\e}}\frac{e^{-i\mbf{u} \cdot \mbf{\xi}}}{-2iu_1}\left(1-\frac{|u_1|}{T}\right)\left(1-\frac{|u_2|}{T}\right)(\phi(\mbf{u})-\chi(\mbf{u}))k_{\e}(\mbf{u}) du_1 du_2, \\
J_{\e,2}(\mbf{\xi}) &:= \frac{T^2}{2\pi} \int_{\mc{S}_{\e}}\frac{e^{-i\mbf{u} \cdot \mbf{\xi}}}{2iu_2}\left(1-\frac{|u_1|}{T}\right)\left(1-\frac{|u_2|}{T}\right)(\phi(\mbf{u})-\chi(\mbf{u}))k_{\e}(\mbf{u}) du_1 du_2.
\end{align*}
Note that the integrand is differentiable almost everywhere. Differentiating under the integral sign gives
\begin{equation*}
h_{T,\e}(\mbf{\xi}) = \frac{\partial J_{\e,1}}{\partial \xi_1}(\mbf{\xi})- \frac{\partial J_{\e,2}}{\partial \xi_2}(\mbf{\xi}).
\end{equation*}
Fix $\tau_1,\tau_2 > 0$, set $A(\mbf{\tau},\mbf{z}) := [-\tau_1,z_1] \times [-\tau_2,z_2]$, and let $R = R_{\tau_1,\tau_2}$ denote the counterclockwise oriented rectangle with corners $(-\tau_1,-\tau_2)$, $(z_1,-\tau_2)$, $(z_1,z_2)$ and $(-\tau_1,z_2)$. By Green's theorem,
\begin{align*}
\int_{A(\mbf{\tau},\mbf{z})} h_{T,\e}(\mbf{\xi})d\xi_1d\xi_2 &= \int_R \left(J_{\e,2}(\mbf{\xi}) d\xi_1 + J_{\e,1}(\mbf{\xi}) d\xi_2\right) \\
&= \int_{-\tau_1}^{z_1} \left(J_{\e,2}(\xi_1,-\tau_2) - J_{\e,2}(\xi_1,z_2)\right)d\xi_1 + \int_{-\tau_2}^{z_2}\left(J_{\e,1}(-\tau_1,\xi_2) - J_{\e,1}(z_1,\xi_2)\right)d\xi_2 \\
&= \frac{T^2}{2\pi} \int_{[-T,T]^2} du_1du_2\left(1-\frac{|u_1|}{T}\right)\left(1-\frac{|u_2|}{T}\right) \left(\phi(\mbf{u})-\chi(\mbf{u})\right) k_{\e}(\mbf{u})\\
&\cdot \left(\left(\frac{e^{-iu_2z_2} - e^{iu_2\tau_2}}{2i} \right)\int_{-\tau_1}^{z_1} \frac{e^{-iu_1\xi_1}}{u_2}d\xi_1 - \left(\frac{e^{iu_1\tau_1} - e^{-iu_1z_1}}{2i}\right)\int_{-\tau_2}^{z_2} \frac{e^{-iu_2\xi_2}}{u_1} d\xi_2\right).
\end{align*}
Evaluating the integrals and bounding trivially then yields
\begin{align*}
&\left|\int_{-\tau_1}^{z_1}\int_{-\tau_2}^{z_2} h_{T,\e}(\mbf{\xi})d\xi_1d\xi_2\right| \\
&\leq \frac{T^2}{2\pi} du_1du_2\int_{[-T,T]^2} du_1du_2\frac{|\phi(\mbf{u})-\chi(\mbf{u})|}{|u_1u_2|}\\
&\cdot\left|\left(\frac{e^{-iu_2z_2} - e^{iu_2\tau_2}}{2} \right)\left(e^{-iu_1z_1}-e^{iu_1\tau_1}\right) -\left(\frac{e^{iu_1\tau_1} - e^{-iu_1z_1}}{2}\right)\left(e^{-iu_2z_2} - e^{iu_2\tau_2}\right)\right| \\
&\ll \int_{[-T,T]^2} \frac{|\phi(\mbf{u})-\chi(\mbf{u})|}{|u_1u_2|}\left|\left(e^{-iu_1(z_1+\tau_1)} - 1\right)\left(e^{-iu_2(z_2+\tau_2)} - 1\right)\right|du_1du_2.
\end{align*}
Note that this holds for all $\e > 0$, so again by the Dominated Convergence theorem, we also have
\begin{equation*}
\left|\int_{-\tau_1}^{z_1}\int_{-\tau_2}^{z_2} h_{T}(\mbf{\xi})d\xi_1d\xi_2\right| \ll \int_{[-T,T]^2} \frac{|\phi(\mbf{u})-\chi(\mbf{u})|}{|u_1u_2|}\left|\left(e^{-iu_1(z_1+\tau_1)} - 1\right)\left(e^{-iu_1(z_1+\tau_1)} - 1\right)\right|du_1du_2.
\end{equation*}
Suppose that $|\tau_1|,|\tau_2| \leq T$ and fix $\mbf{u} \in [-T,T]^2 \bk \mc{R}_T$, and suppose without loss of generality that $|u_1| \leq T^{-3}$ (the alternative case being similarly argued). As $\left|e^{-iu_1(z_1 + \tau_1)} - 1\right| \leq 2|u_1||z_1+\tau_1|$, 
\begin{align*}
\int_{|u_1| \leq T^{-3}} \int_{|u_2| \in [-T^{-3},T]} &\frac{|\phi(\mbf{u})-\chi(\mbf{u})|}{|u_1u_2|}\left|\left(e^{-iu_1(z_1+\tau_1)} - 1\right)\left(e^{-iu_1(z_1+\tau_1)} - 1\right)\right|du_1du_2 \\
&\leq |z_1 + \tau_1|T^{-2} \ll T^{-1}.
\end{align*}
Similarly, if both $|u_1|,|u_2| \leq T^{-3}$, we have
\begin{equation*}
\int_{[-T^{-3},T^{-3}]} \frac{|\phi(\mbf{u})-\chi(\mbf{u})|}{|u_1u_2|}\left|\left(e^{-iu_1(z_1+\tau_1)} - 1\right)\left(e^{-iu_1(z_1+\tau_1)} - 1\right)\right|du_1du_2 \leq |z_1 +\tau_1||z_2 + \tau_2|T^{-6} \ll T^{-1}.
\end{equation*}
Thus, we have
\begin{equation*}
\int_{[-T,T]^2 \bk \mc{R}_T} \frac{|\phi(\mbf{u})-\chi(\mbf{u})|}{|u_1u_2|}\left|\left(e^{-iu_1(z_1+\tau_1)} - 1\right)\left(e^{-iu_1(z_1+\tau_1)} - 1\right)\right|du_1du_2 \ll T^{-1},
\end{equation*}
and hence, choosing $\mbf{\tau} = (T,T)$,
\begin{align*}
|\Delta \ast \mu_T(A((T,T), \mbf{z}))| &\ll \int_{\mc{R}_T} \frac{|\phi(\mbf{u})-\chi(\mbf{u})|}{|u_1u_2|}\left|\left(e^{-iu_1(z_1+\tau_1)} - 1\right)\left(e^{-iu_1(z_1+\tau_1)} - 1\right)\right|du_1du_2 + T^{-1} \\
&\ll \int_{\mc{R}_T} \frac{|\phi(\mbf{u})-\chi(\mbf{u})|}{|u_1u_2|} du_1du_2 + T^{-1}.
\end{align*}
Since the convolution operator $\Delta \mapsto \Delta \ast \mu_T$ has norm $\|W_T\|_{L^1(\mb{R}^2)} \leq T^2$ on the space of probability measures, \eqref{TAILS} implies that
\begin{align*}
|\Delta \ast \mu_T(B(\mbf{z}) \bk A((T,T),\mbf{z}))| &\leq T^2|\Delta(B(\mbf{z}) \bk A((T,T),\mbf{z}))| \\
&= T^2|G(\mbf{z}) - H(\mbf{z}) - \int_{-T}^{z_1}\int_{-T}^{z_2} d(G-H)(\mbf{u})| \ll T.
\end{align*}
Combining this with \eqref{LEFT} and dividing by $T^2$ proves the claim.
%Since this holds for every $\tau_1,\tau_2 > 0$, taking them to infinity yields
%\begin{equation*}
%|\Delta \ast \mu_T(\mbf{z})| \ll \int_{[-T,T]^2} \frac{|\phi(\mbf{u})-\chi(\mbf{u})|}{|u_1u_2|}\left|\left(e^{-iu_1(z_1+\tau_1)} - 1\right)\left(e^{-iu_1(z_1+\tau_1)} - 1\right)\right|du_1du_2.
%\end{equation*}
\end{proof}
We next show that Lemma \ref{ESSEEN} is applicable to the case that $H = F_{x,S,a}$ and $G = \Phi_{(2)}$ when $S$ is siftable with respect to $a$.
\begin{lem} \label{22MOMENTS}
Let $a \in \mb{N}$ and let $S$ be siftable with respect to $a$. Then for $\e_1,\e_2 \in \{0,1\}$,
\begin{equation*}
\sum_{n \leq x \atop (n,a) = 1} 1_S(n)1_S(n+a) \left(\frac{\omega(n)-\lambda(x)}{\sqrt{\lambda(x)}}\right)^{2\e_1}\left(\frac{\omega(n+a)-\lambda(x)}{\sqrt{\lambda(x)}}\right)^{2\e_2} = E(x;a)\left(1+O\left(\frac{1}{\sqrt{\lambda(x)}}\right)\right).
\end{equation*}
\end{lem}
\begin{proof}
We will only prove the case that $\e_1 = \e_2 = 1$, as the remaining cases are either trivial or follow from a simpler version of the argument we present here. \\
We follow the method of Granville and Soundararajan \cite{GrS}. For each prime $p \leq x$, define $g_p(n) := 1-f(p)$ when $p | n$, and $g_p(n) = -f(p)$ when $p \nmid n$. Let $z = x^{\frac{\theta}{4}}$. For any $m \leq x$,
\begin{equation*}
\omega(m)-\lambda(x) = \sum_{p|m} \left(1-f(p)\right) - \sum_{p \leq x \atop p\nmid m} f(p) = \sum_{p \leq x} g_p(m) = \sum_{p \leq y} g_p(m) + O\left(1\right).
\end{equation*}
Moreover, for $m = \prod_j p_j^{\alpha_j}$ put $g_m(n) := \prod_j g_{p_j}(n)^{\alpha_j}$. Then
\begin{align*}
\mu_{2,2}(x) &:= \sum_{n \leq x \atop (n,a) = 1} 1_S(n)1_S(n+a)\left(\omega(n)-\lambda(x)\right)^2\left(\omega(n+a)-\lambda(x)\right)^2 \\
&= \sum_{p_1,p_2,q_1,q_2 \leq z} \sum_{n\leq x \atop (n,a) = 1} 1_S(n) 1_S(n+a) g_{p_1p_2}(n)g_{q_1q_2}(n+a) \\
&+ O\left(\sum_{n \leq x \atop (n,a) = 1} 1_S(n)1_S(n+a) \left|\sum_{p \leq z} g_p(n)\right|\left|\sum_{q \leq z} g_q(n+a)\right|\left(\left|\sum_{p \leq z} g_p(n)\right|+\left|\sum_{q \leq z} g_q(n+a)\right|\right)\right).
\end{align*}
Let $M_{2,2}(x;p_1,p_2,q_1,q_2)$ be the first term above. Since $g_{p_1p_2}(n) = g_{p_1p_2}((n,p_1p_2))$, we have
\begin{align*}
&M_{2,2}(x;p_1,p_2,q_1,q_2) = \sum_{d|p_1p_2}\sum_{e|q_1q_2}g_{p_1p_2}(d)g_{q_1q_2}(e) \sum_{n \leq x \atop (n,p_1p_2) = d, (n+a,q_1q_2) = e} 1_S(n)1_S(n+a)1_{(n,a)= 1} \\
&= \sum_{d | p_1p_2}\sum_{e|q_1q_2} g_{p_1p_2}(d)g_{q_1q_2}(e) \sum_{d'|p_1p_2/d}\mu(d')\sum_{e' | q_1q_2/e}\mu(e') \sum_{n \leq x \atop n \equiv 0 (dd'), n + a \equiv 0 (ee')} 1_S(n)1_S(n+a) 1_{(n,a) = 1}.
\end{align*}
Note that $(dd',ee') = 1$ since $(n,a) = 1$ for each $n$ in the support of the sum. Since $S$ is siftable with respect to $a$, we have
\begin{align*}
M_{2,2}(x;p_1,p_2,q_1,q_2) &= \sum_{d | p_1p_2}\sum_{e|q_1q_2} g_{p_1p_2}(d)g_{q_1q_2}(e)\\
&\cdot\sum_{d'|p_1p_2/d}\mu(d')\sum_{e' | q_1q_2/e}\mu(e')f(dd'ee')\left(E(x;a) + R(x;dd',ee')\right) \\
&= E(x;a)\left(\sum_{dd'| p_1p_2} \mu(d')g_{p_1p_2}(d)f(dd')\right)\left(\sum_{ee'| p_1p_2} \mu(e')g_{q_1q_2}(e)f(ee')\right) \\
&+ \sum_{dd'|p_1p_2}\sum_{ee'|q_1q_2} \mu(d'e')g_{p_1p_2}(d)g_{q_1q_2}(e) f(dd'ee')R(x;dd',ee').
\end{align*}
Since the terms depending on $p_1p_2$ are identical (up to relabelling) to those depending on $q_1q_2$, we shall only write out the treatment of the former terms. We first consider the case that $p_1 \neq p_2$. In that case, $(d,d') = 1$ necessarily, and thus
\begin{align*}
&\sum_{dd'|p_1p_2}\mu(d')g_{p_1p_2}(d)f(dd') =  \sum_{d|p_1p_2} g_{p_1p_2}(d)f(d)\prod_{p|p_1p_2/d}\left(1-f(p)\right) \\
&= \left(\prod_{p|p_1p_2} \left(1-f(p)\right)\right) \sum_{d|p_1p_2}g_{p_1p_2}(d)\prod_{p|d} f(p)\left(1-f(p)\right)^{-1} \\
&= \left(\prod_{p|p_1p_2} \left(1-f(p)\right)\right)\left(f(p_1)f(p_2) -f(p_1)f(p_2) - f(p_1)f(p_2) + f(p_1)f(p_2)\right)= 0.
\end{align*}
It follows that the off-diagonal terms contribute
\begin{align*}
&\sum_{p_1,p_2,q_1,q_2 \leq y \atop p_1 \neq p_2 \text{ or } q_1 \neq q_2} M_{2,2}(x;p_1,p_2,q_1,q_2) \\
&= \sum_{p_1,p_2,q_1,q_2 \leq z \atop p_1 \neq p_2 \text{ or } q_1 \neq q_2} \sum_{dd'|p_1p_2}\sum_{ee'|q_1q_2} \mu(d'e')g_{p_1p_2}(d)g_{q_1q_2}(e) f(dd'ee')R(x;dd',ee') =: \mc{R}_1.
\end{align*}
Conversely, if $p_1 = p_2$ then a simple computation 
%using the condition $f(p^2) \leq Af(p)$ for some $A \geq 1$, 
shows that
\begin{align*}
\sum_{dd'|p_1p_2}\mu(d')g_{p_1p_2}(d)f(dd') &= f(p_1)^2 + f(p_1)(1-f(p_1))^2 - f(p_1)^3\\
&= f(p_1)+O(f(p_1)^2).
\end{align*}
The diagonal terms thus contribute
\begin{align*}
\sum_{p,q \leq z} M_{2,2}(x;p,p,q,q) &= E(x;a) \left(\sum_{p \leq z} f(p) + O(1)\right)^2 \\
&+ \sum_{p,q \leq z} \sum_{dd'|p^2}\sum_{ee'|q^2} \mu(d'e')g_{p^2}(d)g_{q^2}(e) f(dd'ee')R(x;dd',ee') \\
&= E(x;a) \left(\lambda(x) + O(1)\right)^2 + \mc{R}_2.
\end{align*}
Now, for any $m \in \mb{N}$, $|g_m| \leq 1$, so that since $\beta \ra \infty$,
\begin{equation*}
\mc{R}_1 \ll \sum_{p_1,p_2,q_1,q_2 \leq z} \max_{dd'|p_1p_2 \atop ee'|q_1q_2} |R(x;dd',ee')| \leq \sum_{m_1m_2 \leq z^4} |R(x;m_1,m_2)| \ll x(\log x)^{-3}.
\end{equation*}
We can estimate $\mc{R}_2$ the same way. It follows that
\begin{align*}
\mu_{2,2}(x) &= E(x;a)\lambda(x)^2\left(1+O\left(\frac{1}{\lambda(x)}\right)\right) \\
&+ O\left(\sum_{n \leq x \atop (n,a) = 1} 1_S(n)1_S(n+a) \left|\sum_{p \leq z} g_p(n)\right|\left|\sum_{q \leq z} g_q(n+a)\right|\left(\left|\sum_{p \leq z} g_p(n)\right|+\left|\sum_{q \leq z} g_q(n+a)\right|\right)\right).
\end{align*}
We will estimate $$\sum_{n \leq x \atop (n,a) = 1}1_S(n)1_S(n+a)\left|\sum_{p \leq z} g_p(n)\right|^2\left|\sum_{q \leq z} g_q(n+a)\right|,$$ the remaining term being estimated in exactly the same way. By the Cauchy-Schwarz inequality,
\begin{align*}
&\sum_{n \leq x \atop (n,a) = 1}1_S(n)1_S(n+a)\left|\sum_{p \leq z} g_p(n)\right|^2\left|\sum_{q \leq z} g_q(n+a)\right| \\
&\leq \left(\sum_{n \leq x \atop (n,a) = 1} 1_S(n)1_S(n+a) \left(\sum_{p \leq z} g_p(n)\right)^2\left(\sum_{q \leq z} g_q(n+a)\right)^2\right)^{\frac{1}{2}} \\
&\cdot\left(\sum_{n \leq x \atop (n,a) = 1} 1_S(n)1_S(n+a)\left(\sum_{p \leq z} g_p(n)\right)^2\right)^{\frac{1}{2}}.
\end{align*}
The first sum here is precisely the $(2,2)$-moment we just estimated, and its square root is thus $\ll \lambda(x)\sqrt{E(x;a)}$. The second is the $(2,0)$-moment, i.e., where $\e_1 = 1$ and $\e_2 = 0$, and its square root is $\ll \sqrt{\lambda(x)E(x;a)}$. Hence,
\begin{equation*}
\mu_{2,2}(x) = E(x;a) \lambda(x)^2\left(1+O\left(\frac{1}{\sqrt{\lambda(x)}}\right)\right).
\end{equation*}
Dividing both sides by $\lambda(x)^2$ now suffices to prove the claim.
\end{proof}
\begin{cor} \label{APPLY}
If $a \in \mb{N}$ and $S$ is siftable with respect to $a$ then for each fixed $\mbf{z}$ and $T \geq 1$,
\begin{equation*}
\left|F_{x,S,a}(\mbf{z})-\Phi_{(2)}(\mbf{z}) - \int_{-T}^{z_1}\int_{-T}^{z_2} d(F_{x,S,a}-\Phi_{(2)})(\mbf{u})\right| \ll T^{-1}.
\end{equation*}
\end{cor}
\begin{proof}
Observe that
\begin{equation*}
\left|\Phi_{(2)}(\mbf{z}) - \int_{-T}^{z_1}\int_{-T}^{z_2} d\Phi_{(2)}(\mbf{w})\right| \ll \int_{-\infty}^{-T} \int_{\mb{R}} e^{-\frac{1}{2}(u^2+v^2)} du dv \ll e^{-\frac{1}{4}T^2} \ll T^{-1}.
\end{equation*}
Moreover, suppose that $n \leq x$ is counted by $F_{x,S,a}(\mbf{z}) - \int_{-T}^{z_1}\int_{-T}^{z_2} dF_{x,S,a}(\mbf{w})$. Then $\min\{\tilde{\omega}_x(n),\tilde{\omega}_x(n+a)\} < -T$. It follows by Lemma \ref{22MOMENTS} that
\begin{align*}
&\left|F_{x,S,a}(\mbf{z}) - \int_{-T}^{z_1}\int_{-T}^{z_2} dF_{x,S,a}(\mbf{w})\right| \\
&\ll \sum_{\e \in \{0,1\}}E(x;a)^{-1}|\{n \leq x : n,n+a \in S, (n,a) = 1, \left|\frac{\omega(n+\e a)-\lambda(x)}{\sqrt{\lambda(x)}}\right| > T\}| \\
&\ll T^{-2}E(x;a)^{-1} \sum_{n \leq x \atop (a,n) = 1} 1_S(n)1_S(n+a) \left(1+\left(\frac{\omega(n)-\lambda(x)}{\sqrt{\lambda(x)}}\right)^2\right)\left(1+\left(\frac{\omega(n+a)-\lambda(x)}{\sqrt{\lambda(x)}}\right)^2\right) \\
&\ll T^{-2}.
\end{align*}
This more than suffices to prove the claim.
\end{proof}
%\begin{proof}
%%Notice that 
%%\begin{align*}
%%F_{x,S,a}(\mbf{z})-\Phi_{(2)}(\mbf{z}) &= \mb{P}(\tilde{\omega}_x> z_1,\tilde{\omega}_x \circ \sg^a \leq z_2) - \left(1-\Phi(z_1)\right)\Phi(z_2) + \mb{P}(\tilde{\omega}_x \leq  z_1,\tilde{\omega}_x \circ \sg^a > z_2) - \Phi(z_1)\left(1-\Phi(z_2)\right) \\
%%&= \
%It obviously suffices to show that
%\begin{equation*}
%\mb{P}\left(|\tilde{\omega}_{x,S}| > z_1, |\tilde{\omega}_{x,S} \circ \sg^a| > z_2\right) \ll (z_1z_2)^{-2}.
%\end{equation*}
%But this is trivial from Chebyshev's inequality, in light of Lemma \ref{22MOMENTS}.
%\end{proof}
The following simple lemma allows us to put into effect the Poisson approximation results of Section 3.
\begin{lem} \label{TRIV} 
Let $W,Z$ be random vectors on a common probability space taking values in $\mb{R}^n$ with respective laws $\mu$ and $\nu$. Then
\begin{equation*}
\sup_{\mbf{u} \in \mb{R}^n} \left|\mb{E}\left(e^{i\mbf{u} \cdot W}\right)-\mb{E}\left(e^{i\mbf{u} \cdot Z}\right)\right| \ll d_{TV}(\mu,\nu).
\end{equation*}
\end{lem}
\begin{proof}
Fix $\mbf{u} \in \mb{R}^n$. For each $r \in [0,2\pi)$ let $A(r) := \{\mbf{t} \in \mb{R}^n : \mbf{t} \cdot \mbf{u} \equiv r (2\pi)\}$, which is clearly Borel measurable. Then
\begin{align*}
\left|\mb{E}\left(e^{i\mbf{u} \cdot W}\right)-\mb{E}\left(e^{i\mbf{u} \cdot Z}\right)\right| &= \left|\int_0^{2\pi} dre^{ir}\int_{A(r)} d(\mu-\nu)(\mbf{t})\right| \leq 2\pi \sup_{r \in [0,2\pi)} |\mu(A(r))-\nu(A(r))| \\
&\ll d_{TV}(\mu,\nu),
\end{align*}
as claimed.
\end{proof}
The following is a trivial calculation.
\begin{lem} \label{COMPUTE}
Let $X$ be a Poisson random variable with parameter $\lambda$, and let $\tilde{X} := (X-\lambda)/\sqrt{\lambda}$. Let $\pi$ denote the characteristic function of $\tilde{X}$. Then
\begin{equation*}
|\pi(w)| = \exp\left(-\lambda(1-\cos(w/\sqrt{\lambda}))\right).
\end{equation*}
\end{lem}
We thus have the following.
\begin{prop}\label{PHIEST}
Let $a \in \mb{N}$ and let $S$ be a siftable, regular set (with respect to $a$). Let $|u|,|v| \leq 2\pi\sqrt{\log_2 x}$, and let $w := \max\{|u|,|v|\}$.  Then
\begin{align}
|\phi_{x,y}(u,v) - \chi(u)\chi(v)| &\ll \frac{1}{\log_2 y} + R(x;y) \label{B1} \\
|\phi_x(u,v)-\chi(u)\chi(v)| &\ll \frac{w\log_2(\mf{L}/w)}{(\log_2 x)^{\frac{1}{2}}}. \label{B4}
%|\phi_x(u,v) - \pi_z(u)\pi_z(v)| \ll \frac{w\log(\sqrt{\log_2 x}/w)}{\sqrt{\log_2 x}},
\end{align}
%as well as
%\begin{equation}\label{PHIXZ}
%|\phi_{x,z}(u,v)| \ll \frac{1}{\log_2 z} + R(x;z) + \exp\left(-\lambda(2-\cos(u/\sqrt{\lambda}) -\cos(v/\sqrt{\lambda}))\right),
%\end{equation}
%and also
%\begin{align}
%|\phi_{y}(u,v)| &\ll \frac{1}{\log_2 z} + R(x;z) + \exp\left(-\frac{1}{3}(u^2+v^2)\right) \label{B2}\\
%|\phi_x(u,v)| &\ll \frac{w\log(\sqrt{\log_2 x}/w)}{\sqrt{\log_2 x}}+ \exp\left(-\frac{1}{3}(u^2 + v^2)\right). \label{B3}
%%+ \frac{\max \{|u|,|v|\}\log \beta}{\sqrt{\log_2 x}}.
%\end{align}
\end{prop}
\begin{proof}
Let $Z := (Z_1,Z_2) \sim \text{Poi}(\lambda(y))^2$ be the independent random variables from Lemma \ref{TVPOISSON}. Let $\tilde{Z}_j := (Z_j-\lambda(y))/\sqrt{\lambda(y)}$, $\tilde{Z} := (\tilde{Z}_1,\tilde{Z}_2)$ and let $\pi_y$ denote the characteristic function of $\tilde{Z}_j$. Since $S$ is regular, we note that $\lambda(y) \gg \log_2 y$. Note that by independence, we have $\mb{E}\left(e^{i(u,v) \cdot \tilde{Z}}\right) = \pi_y(u)\pi_y(v)$. By Lemmata \ref{TVPOISSON} and \ref{TRIV}, we have
\begin{equation} \label{POICOMP}
|\phi_{x,y}(u,v)-\pi_y(u)\pi_y(v)| \ll d_{TV}(\nu_y,\mc{L}(\tilde{Z})) \ll \frac{1}{\log_2 y} + R(x;y).
\end{equation}
%By the triangle inequality and Lemma \ref{COMPUTE}, it follows that
%\begin{align*}
%|\phi_{y}(u,v)| &\leq |\phi_{y}(u,v)-\pi_y(u)\pi_y(v)| + |\pi_y(u)||\pi_y(v)| \\
%&\ll \frac{1}{\log_2 y} + R(x;y) + \exp\left(-\lambda(2-\cos(u/\sqrt{\lambda}) -\cos(v/\sqrt{\lambda}))\right) \\
%&\leq \frac{1}{\log_2 y} + R(x;y) + \exp\left(-\frac{1}{3}(u^2 + v^2)\right),
%\end{align*}
%by using the simple lower bound $1-\cos t \geq t^2/3$. This gives \eqref{B2}.\\
%Invoking Lemma \ref{TVPOISSON} then establishes the first assertion. \\
%thus establishing \eqref{B3}. \\
We now compare $\pi_y(u)\pi_y(v)$ to $\chi(u)\chi(v)$. We note that when $|v^4/\lambda(y)| < 1$ we have
\begin{equation*}
|\pi_y(v)-\chi(v)| = \left|\exp\left(-\lambda(y)(1-\cos(v/\sqrt{\lambda(y)}))\right)-e^{-\frac{1}{2}v^2}\right| \ll e^{-\frac{1}{2}v^2} v^4\lambda(y)^{-1}.
\end{equation*}
On the other hand, when $|v^4/\lambda(y)| > 1$ then
\begin{equation*}
|\pi_y(v) - \chi(v)| \leq |\pi_y(v)| + e^{-\frac{1}{2}v^2} \ll e^{-\frac{1}{3}\sqrt{\lambda(y)}}.
\end{equation*}
It follows that
\begin{align*}
|\pi_y(u)\pi_y(v) - \chi(u)\chi(v)| &\ll \max\{|\pi_y(v)-\chi(v)|,|\pi_y(u)-\chi(u)|\} \\
&\ll e^{-\frac{1}{2}w^2}w^4\lambda(y)^{-1} 1_{|w| < \lambda(y)^{1/4}} + e^{-\frac{1}{3}\sqrt{\lambda(y)}}1_{|w| > \lambda(y)^{1/4}}.
\end{align*} 
This, coupled with \eqref{POICOMP}, implies \eqref{B1}. \\
Now, we consider the difference between $\phi_x$ and $\phi_{x,y}$. We have
%We first observe that by Lemma \ref{22MOMENTS},
%\begin{equation*}
%|\{n \leq x : n,n+a \in S, (n,a) = 1, |\omega(n)-\lambda(x)| > Z\sqrt{\lambda(x)}\}| \ll E(x;a)Z^{-2},
%\end{equation*}
%so choosing $Z = \lambda(x)^{1/6}$, we see that
\begin{align*}
|\phi_x(u,v)-\phi_{x,y}(u,v)| &\leq E(x;a)^{-1}\sum_{n \leq x \atop (n,a) = 1} 1_S(n)1_S(n+a)\left|e^{iu(\tilde{\omega}_x(n)-\tilde{\omega}_y(n))}e^{iv(\tilde{\omega}_x(n+a)-\tilde{\omega}_y(n+a))}-1\right| \\
&\ll E(x;a)^{-1}\sum_{n \leq x \atop (n,a) = 1} 1_S(n)1_S(n+a) \\
&\cdot \max\{|u||\tilde{\omega}_x(n)-\tilde{\omega}_y(n)|,|v||\tilde{\omega}_x(n+a)-\tilde{\omega}_y(n+a)|\} \\
&\ll \max\{|u|,|v|\} E(x;a)^{-1}\sum_{n \leq x\atop (n,a) = 1} 1_S(n)1_S(n+a) \\
&\cdot \left(\frac{|\omega(n)-\omega_y(n)|-|\lambda(x)-\lambda(y)|}{\sqrt{\lambda(y)}} + \left(1-\sqrt{\frac{\lambda(y)}{\lambda(x)}}\right)|\omega(n)-\lambda(x)|\right) \\
&\ll  w\left(\frac{\beta}{\sqrt{\log_2 x}} + \left(\frac{\log \beta}{\log_2 x}\right)E(x;a)^{-1}\sum_{n \leq x\atop (n,a) = 1} 1_S(n)1_S(n+a)|\omega(n)-\lambda(x)|\right).
\end{align*}
By the Cauchy-Schwarz inequality and Lemma \ref{22MOMENTS}, we have
\begin{equation*}
E(x;a)^{-1}\sum_{n \leq x} 1_S(n)1_S(n+a)1_{(n,a) = 1}|\omega(n)-\lambda(x)| \ll \sqrt{\lambda(x)},
\end{equation*}
so that $|\phi_x(u,v) - \phi_{x,y}(u,v)| \ll w \beta (\log_2 x)^{-\frac{1}{2}}$. 
%This, coupled with \eqref{B2} thus gives
%\begin{equation*}
%|\phi_x(u,v)| \ll \max\{|u|,|v|\} \frac{\log \beta}{\mf{L}} + \frac{1}{\log_2 y} + e^{-\frac{1}{6}\beta \log \beta}.
%\end{equation*}
Choose $y$ such that $\beta \log \beta = 21\log\left(\mf{L}/w\right)$. Then, upon applying \eqref{B1},
\begin{align*}
|\phi_x(u,v)-\chi(u)\chi(v)| &\leq |\phi_x(u,v)-\phi_{x,y}(u,v)| + |\phi_{x,y}(u,v)-\chi(u)\chi(v)|\\
&\ll \frac{w\log_2(\mf{L}/w)}{(\log_2 x)^{\frac{1}{2}}} + \frac{1}{\log(\log x \log(\mf{L}/w))} + e^{-\frac{1}{2}w^2} w^4(\log(\log x \log(\mf{L}/w)))^{-1},
\end{align*}
which completes the proof of \eqref{B4}.
\end{proof}
\begin{proof}[Proof of Theorem \ref{EKTHM}]
Let $T := \left(\log_2 x\right)^{\frac{1}{4}}(\log_3 x)^{-1}$. By Lemma \ref{ESSEEN}, we have
\begin{equation*}
\|F_x-\Phi_{(2)}\|_{L^{\infty}(\mb{R}^2)} \ll \int_{\mc{R}_T} \frac{|\phi_x(u,v)-\chi(u)\chi(v)|}{|uv|} dudv + T^{-1},
\end{equation*}
where we recall that $\mc{R}_T := \{\mbf{u} \in [-T,T]^2 : |u_1|,|u_2| \geq T^{-3}\}$. From \eqref{B4}, we have $|\phi_x(u,v) - \chi(u)\chi(v)| \ll w\frac{\log_3 x}{\sqrt{\log_2 x}}$ for $w := \max\{|u|,|v|\}$. Thus, 
%We split the integral into the regions $|u|,|v| \in [T^{-3},1]$ and $\max\{|u|,|v|\} \in [1,T]$. In the first region we use \eqref{B4} and write
\begin{align*}
&\int_{\mc{R}_T} \frac{|\phi_x(u,v)-\chi(u)\chi(v)|}{|uv|} dudv \leq 2\int_{|u| \in [T^{-3},T]} \int_{|v| \leq |u|} \frac{|\phi_x(u,v)-\chi(u)\chi(v)|}{|uv|} dudv \\
&\ll \frac{\log_3 x}{\sqrt{\log_2 x}} \int_{|u| \in [T^{-3},T]} du \int_{|v| \in [T^{-3},T]} \frac{dv}{v} \ll \frac{\log_3 x}{\sqrt{\log_2 x}}(T \log T) = (\log_3 x)(\log_2 x)^{-\frac{1}{4}}.
\end{align*}
%upon making the change of variables $u \mapsto 1/u$. \\
%%\ref{UNIF}, giving
%%\begin{equation*}
%%|\phi_x(u,v)-\chi(u)\chi(v)| \ll \log_2^{-\frac{1}{4}} x
%%\end{equation*}
%%for all $|u|,|v| \in [T^{-2},1]$ and $x$ sufficiently large. Selecting $T := \log_2^{\frac{1}{4}} x$ gives
%%\begin{equation*}
%%\int_{|u|,|v| \in [T^{-2},1]} \frac{|\phi_x(u,v)-\chi(u)\chi(v)|}{uv} dudv \ll \log^2(T^2) \log_2^{-\frac{1}{4}} x \ll \log_3^2 x\left(\log_2 x\right)^{-\frac{1}{4}}.
%%\end{equation*}
%Call the remaining region $\mc{R}_T'$.  for $(u,v) \in \mc{R}_T'$. As such, 
%\begin{align*}
%\int_{\mc{R}} \frac{|\phi_x(u,v)-\chi(u)\chi(v)|}{uv} &\ll \left(\int_{|u| \in [T^{-3},1]} + \int_{|u| \in [1,T]}\right)du \int_{|v| \in [1,T]}dv \frac{|\phi_x(u,v)-\chi(u)\chi(v)|}{uv} \\
%&\ll \log_3 x\left(\log_2 x\right)^{-\frac{1}{3}} T\log T + \frac{\log_3 x}{(\log_2 x)^{\frac{1}{3}}}\ll \log_3 x\left(\log_2 x\right)^{-\frac{1}{6}}.
%\end{align*}
Thus, we get $\|F_x-\Phi^2\|_{L^{\infty}(\mb{R}^2)} \ll (\log_3 x)\left(\log_2 x\right)^{-\frac{1}{4}}$, and the proof of Theorem \ref{EKTHM} is complete in the case that $1_S$ is not multiplicative. \\
When $1_S$ is multiplicative and $S$ is siftable for each $\alpha | a$ then as $1_S(n)1_S(n+a) = 1_S(m)1_S(m+\gamma)$, where $\alpha \gamma = a$ and $m\alpha = n$,
\begin{align*}
\phi^{\ast}_{x}(u,v) &= E^{\ast}(x;a)^{-1}\sum_{\alpha \gamma = a} e^{i(u+v)\lambda(x)^{-\frac{1}{2}} \omega(\alpha)}  \sum_{m \leq x/\alpha \atop (m,\gamma) = 1} 1_S(m)1_S(m+\gamma) e^{i(u\tilde{\omega}_x(m) + v\tilde{\omega}_x(m+\gamma))} \\
&= \sum_{\alpha \gamma = a} e^{i(u+v)\lambda(x)^{-\frac{1}{2}}\omega(\alpha)} \frac{E(x/\alpha;\gamma)}{E^{\ast}(x;a)} \phi_{x/\alpha,S,\gamma}(u,v).
\end{align*}
The result now follows by applying the arguments for the previous case to the average pointwise distance $|\phi_{x/\alpha,S,\gamma}(u,v) - \chi(u,v)|$ for each $\alpha | a$.
%If $1_S$ is multiplicative then we note that
%\begin{equation*}
%\sum_{n \leq x} 1_S(n)1_S(n+a) e^{i(u\omega(n) + v\omega(n+a))} = \sum_{\alpha\beta = a} e^{i\omega(\alpha)(u+v)}\sum_{m \leq x/\alpha \atop (m,\beta) = 1} 1_S(m)1_S(m+\beta) e^{i(u\omega(m) + v\omega(m+\beta))}
\end{proof}
\section{A Sieve Result for Consecutive Squarefree Integers}
Throughout this section, let $E(x) := x\prod_p\left(1-2p^{-2}\right)$ and for $q \geq 2$, let $E_q(x) := E(x)\prod_{p |q} \left(1-2p^{-2}\right)^{-1}$. It is well-known that $E(x)$ is asymptotically the number of integers $n \leq x$ such that $n$ and $n+1$ are both squarefree (see, for instance, \cite{HB}). We show here that $S$ is siftable and regular with respect to each fixed $a \in \mb{N}$, in the sense of Definition \ref{SIFTABLE}. \\Since property ii) in Definition \ref{SIFTABLE} requires information about the distribution of squarefree integers in certain arithmetic progressions, we begin with the following.
%; see Section 4 for a sense of how this will be applied.
\begin{lem} \label{HB}
Let $q \geq 2$ and let $c$ be a reduced residue class modulo $q$, and suppose that $(q,c)$ is squarefree. Let $a \geq 1$. Then
\begin{equation*}
\sum_{n \leq x \atop n \equiv c (q)} \mu^2(n)\mu^2(n+a) = \frac{E_q(x)}{q}\prod_{p || (q,c(c+a))} \left(1-\frac{1}{p}\right)\prod_{p^2||(q,c+a)} \left(1-\frac{1}{p^2}\right)\prod_{p^2||a} \left(1+\frac{1}{p^2}\right) + O\left(x^{\frac{2}{3}+\e} \tau(q)\right).
\end{equation*}
\end{lem}
The restriction that $(q,c)$ be squarefree makes the result non-trivial, as if $q$ and $c$ shared a common factor $b^2$ say, then if $n \equiv c(q)$ then $b^2 |n$. Thus, the sum here would be 0.
\begin{proof}
For $e_1,e_2 \leq x^{\frac{1}{2}}$ let 
\begin{equation*}
N_a(x;e_1,e_2,c,q) := \left|\{n \leq x : e_1^2|n, e_2^2 |(n+a), n \equiv c(q)\}\right|.
\end{equation*}
Note that $N_a(x;e_1,e_2,c,q)$ is zero as long as $(e_1,e_2)^2 \nmid a$. Let $2 \leq z \leq x^{\frac{1}{2}}$ be a parameter to be chosen. Then
\begin{align*}
&\sum_{n \leq x \atop n \equiv c (q)} \mu^2(n)\mu^2(n+a) = \sum_{e_1,e_2 \leq x^{\frac{1}{2}} \atop (e_1,e_2)^2 |a} \mu(e_1)\mu(e_2) N_a(x;e_1,e_2,c,q) \\
&= \sum_{e_1e_2 \leq z \atop (e_1,e_2)^2 |a} \mu(e_1)\mu(e_2)N_a(x;e_1,e_2,c,q) + \sum_{e_1,e_2 \leq x^{\frac{1}{2}} \atop e_1e_2 > z, (e_1,e_2)^2 |a} \mu(e_1)\mu(e_2)N_a(x;e_1,e_2,c,q) \\
&=: T_1 + T_2.
\end{align*}
Observe that $N_a(x;e_1,e_2,c,q)$ counts simultaneous solutions to the triple of congruences $n \equiv c (q)$, $n \equiv 0 (e_1^2)$ and $n \equiv -a (e_2^2)$. By the Chinese remainder theorem, a solution exists modulo $[q,e_1^2,e_2^2]$ if, and only if, $(q,e_1^2)| c$ and $(q,e_2^2)|(c+a)$ (and $(e_1,e_2)^2 | a$, as we are already assuming). As such, we have
\begin{align*}
T_1 &= x\sum_{e_1e_2 \leq z \atop (e_1,e_2)^2 |a, e_1^2|c,e_2^2|(c+a)} \frac{\mu(e_1)\mu(e_2)}{[q,e_1^2,e_2^2]} + O\left(\sum_{e_1e_2 \leq z}1\right) \\
&= \frac{x}{q}\sum_{e_1e_2 \leq z \atop (e_1,e_2)^2 |a, (q,e_1^2)|c,(q,e_2^2)|(c+a)} \frac{\mu(e_1)\mu(e_2)}{(e_1e_2)^2}(e_1,e_2)^2\left(q,\left(\frac{e_1e_2}{(e_1,e_2)}\right)^2\right) + O\left(\sum_{m \leq z}\tau(m)\right) \\
&= \frac{x}{q}\sum_{\delta^2 | a} \delta^2\sum_{e_1e_2 \leq y \atop (e_1,e_2)=\delta, (q,e_1^2)|c,(q,e_2^2)|(c+a)} \frac{\mu(e_1)\mu(e_2)}{(e_1e_2)^2}(q,e_1^2)(q,(e_2/\delta)^2) + O\left(z \log z\right) \\
&= \frac{x}{q}\sum_{\delta^2 | a} \delta^2\sum_{f_1|(q,c), f_2|(q,c+a)} f_1f_2\sum_{e_1e_2 \leq z \atop (e_1,e_2) = \delta, (e_1^2,q) = f_1, ((e_2/\delta)^2,q) = f_2} \frac{\mu(e_1)\mu(e_2)}{(e_1e_2)^2} + O\left(z \log z\right) \\
&= \frac{x}{q}\sum_{\delta^2 | a} \delta^2\sum_{f_1|(q,c),f_2|(q,c+a)} f_1f_2\sum_{(e_1,e_2) = \delta \atop ((e_1/\delta)^2,q) = f_1, ((e_2/\delta)^2,q) = f_1} \frac{\mu(e_1)\mu(e_2)}{(e_1e_2)^2} + O\left(z \log z + \frac{x}{qz} \tau(q)\right).
\end{align*}
Decompose $(q,c+d) = uv$, where $(u,v) = 1$, $\mu^2(v) = 1$ and $u$ is squarefull. As such, we can split $f_2 = gh$ where $g := (f_2,u)$ and $h := (f_2,v)$. Further, write $\tilde{e}_2 := e_2/\delta$. Note that if $(\tilde{e}_1^2,q) |c$ then $(\tilde{e}_1^2,(q,c)) = (\tilde{e}_1,q)$ (since $(q,c)$ is squarefree), and similarly $(\tilde{e}_2^2,q) = (\tilde{e}_2^2,(q,c+a))$. In particular, $(\tilde{e}_2^2,u) = g$ and $(\tilde{e}_2^2,v) = (e_2,v) = h$. Also, $\mu(\tilde{e}_2) \neq 0$ if, and only if, $g$ is the square of a squarefree integer. Thus, define $G$ implicitly by $G^2 = g$. Then
\begin{align*}
T_1 
%&= \frac{x}{q}\sum_{\delta^2 | d} \delta^2\sum_{G_j^2 | a_j \atop \mu^2(G_j) = 1} \frac{\mu(G_1G_2)(G_1G_2)^2}{(G_1G_2)^4}\sum_{h_j|b_j} \frac{\mu(h_1h_2)h_1h_2}{(h_1h_2)^2} \mathop{\sum_{(e_1,q)= 1}\sum_{((e_2/\delta),q)= 1}}_{(e_1,e_2) = \delta} \frac{\mu(e_1)\mu(e_2)}{(e_1e_2)^2}+ O\left(y \log y + \frac{x}{qy} \tau((q,c(c+d))\right) \\
&= \frac{x}{q}\sum_{\delta^2 | a} \delta^2\sum_{G^2 | u \atop \mu^2(G) = 1} \frac{\mu(G)}{G^2}\sum_{h_1|f_1, h_2|v} \frac{\mu(h_1h_2)}{h_1h_2} \mathop{\sum_{(e_1,q)= 1}\sum_{(\tilde{e}_2,q)= 1}}_{(e_1,e_2) = \delta} \frac{\mu(e_1)\mu(e_2)}{(e_1e_2)^2} + O\left(z \log z + \frac{x}{qz} \tau(q)\right).
\end{align*}
Let $\Sigma_{\delta}$ denote the double sum in $e_1$ and $e_2$ with $(e_1,e_2) = \delta$. Note that if $\mu^2(\delta) = 0$ then $\Sigma_{\delta} = 0$, so assume otherwise. We have
\begin{align*}
\Sigma_{\delta} &= \sum_{(e_1/\delta,q) = 1} \frac{\mu(e_1)}{e_1^2}\sum_{(e_2/\delta,q) = 1, (e_2,e_1) = \delta} \frac{\mu(e_2)}{e_2^2} = \frac{1}{\delta^4}\sum_{(e_1,q) = 1} \frac{\mu(e_1)}{e_1^2}\sum_{(e_2,qe_1) = 1} \frac{\mu(e_2)}{e_2^2} \\
&= \frac{1}{\delta^4}\sum_{(e_1,q) = 1} \frac{\mu(e_1)}{e_1^2} \prod_{p \nmid e_1q} \left(1-\frac{1}{p^2}\right) \\
&= \frac{1}{\delta^4}\prod_{p\nmid q}\left(1-\frac{1}{p^2}\right) \sum_{(e_1,q) = 1} \frac{\mu(e_1)}{e_1^2} \prod_{p|e_1} \left(1-\frac{1}{p^2}\right)^{-1} = \frac{1}{\delta^4}\prod_{p \nmid q}\left(1-\frac{1}{p^2}\right) \prod_{p\nmid q} \left(1-\frac{1}{p^2-1}\right) \\
&= \frac{1}{\delta^4}\prod_{p \nmid q} \left(1-\frac{2}{p^2}\right).
\end{align*}
Inserting this back into our estimate for $T_1$ gives
\begin{align*}
T_1 &= \frac{x}{q}\prod_{p \nmid q} \left(1-\frac{2}{p^2}\right)\sum_{\delta^2 | a \atop \mu^2(\delta) = 1}\frac{1}{\delta^2}\sum_{G^2 | u \atop \mu^2(G) = 1} \frac{\mu(G)}{G^2}\sum_{h_1|f_1, h_2|v} \frac{\mu(h_1h_2)}{h_1h_2} + O\left(z \log z + \frac{x}{qz} \tau(q)\right) \\
&= \frac{x}{q}\prod_{p \nmid q} \left(1-\frac{2}{p^2}\right) \prod_{p || (q,[c,c+a])} \left(1-\frac{1}{p}\right)\prod_{p^2||(q,c+a)} \left(1-\frac{1}{p^2}\right)\prod_{p^2||a}\left(1+\frac{1}{p^2}\right) \\
&+ O\left(z \log z + \frac{x}{qz} \tau(q)\right).
\end{align*}
We next estimate $T_2$. Decomposing the ranges of $e_1$ and $e_2$ dyadically, we have
\begin{align*}
T_2 &\ll \sum_{k,l \leq \frac{\log x}{2 \log 2} \atop k + l \geq \log z} \sum_{2^k < e_1 \leq 2^{k+1} \atop 2^l < e_2 \leq 2^{l+1}} N_a(x;e_1,e_2,c,q) =: \sum_{k,l \leq \frac{\log x}{2 \log 2} \atop k + l \geq \log z} M_{k,l}(x) \ll M_{K,L}(x)\log^2 x 
%\max_{k,l \leq \frac{\log x}{2 \log 2} \atop k+l \geq \log y} \sum_{2^k < e_1 \leq 2^{k+1} \atop 2^l < e_2 \leq 2^{l+1}} N(x;e_1,e_2,c,q),
\end{align*}
where $K$ and $L$ are the respective values of $k$ and $l$ that maximize $M_{k,l}(x)$. Assume without loss of generality that $K \geq L$, the alternative case being similar. Then
\begin{equation*}
M_{K,L}(x) \leq \sum_{2^L < e_2 \leq 2^{L+1}} \sum_{m \leq x/2^{2K}}\sum_{2^K < e_1 \leq 2^{K+1} \atop me_1^2 \equiv -a (e_2^2), me_1^2 \equiv c (q)} 1.
\end{equation*}
The number of solutions in $e_1$ to the simultaneous congruence conditions is at most $\tau(e_2)\tau(q)$, since there are at most 2 solutions to these same congruences modulo each prime $p$ dividing $[e_2,q]$, and each such solution lifts uniquely to a solution mod $p^2$ by Hensel's lemma. Thus, we have
\begin{align*}
M_{K,L}(x) &\leq x2^{-2K}\sum_{l|(q,c+a)} \sum_{2^L < e_2 \leq 2^{L+1} \atop (e_2^2 , q) = l} \left(1+2^{K-2L}lq^{-1}\right) \tau(e_2) \ll xL 2^{L-2K}\left(1+\tau(q)2^{K-2L}q^{-1}\right) \\
&\ll x\log x\left(2^{-K} + \tau(q)2^{-(K+L)}q^{-1} \right) \ll xz^{-\frac{1}{2}} \log x,
\end{align*}
since $K \geq \frac{1}{2}\log z$. As such, we have $T_2 \ll xz^{-\frac{1}{2}} \log^3 x$, and
\begin{align*}
\sum_{n \leq x \atop n \equiv c(q)} \mu^2(n)\mu^2(n+a) &= \frac{x}{q}\prod_{p \nmid q} \left(1-\frac{2}{p^2}\right) \prod_{p || (q,[c,c+a])} \left(1-\frac{1}{p}\right)\prod_{p^2||(q,c+a)} \left(1-\frac{1}{p^2}\right)  \prod_{p^2||a} \left(1+\frac{1}{p^2}\right)\\
&+ O\left(z \log z + \frac{x}{\sqrt{z}} \log^3 x+ \frac{x}{qz} \tau((q,c(c+a)))\right).
\end{align*}
Choosing $z = x^{\frac{2}{3}}$ furnishes a bound stronger than the claim.
\end{proof}
\begin{lem} \label{SIEVESTUFF}
Let $a \geq 1$. Let $q,r \geq 2$ be squarefree, such that $(q,r) |a$. Then
\begin{equation}\label{PLAIN}
\sum_{n \leq x \atop n \equiv 0 (q), n+a \equiv 0 (r)} \mu^2(n)\mu^2(n+a) = x\frac{\phi([q,r])}{[q,r]^2}\prod_{p^2|| a} \left(1+\frac{1}{p^2}\right)\prod_{p \nmid [q,r]} \left(1-\frac{2}{p^2}\right) + O\left(x^{\frac{2}{3}+\e}\tau([q,r])\right).
\end{equation}
Consequently, if $(k,r) = 1$ then 
\begin{equation}\label{NONPLAIN}
\sum_{n \leq x \atop n \equiv 0 (q), n+a \equiv 0 (r)} \mu^2(n)\mu^2(n+a)1_{(n,a) = 1} = E(x;a) \prod_{p|qr} \frac{1}{p}\left(1-\frac{1}{p}\right)\left(1-\frac{2}{p^2}\right)^{-1} + O_a\left(x^{\frac{2}{3}+\e}\tau([q,r])\right),
\end{equation}
where here $$E(x;a) := E(x)\prod_{p^2||a}\left(1+1/p^2\right)\left(1-\frac{1}{p}\left(1-1/p\right)\left(1-2/p^2\right)^{-1}\right).$$
\end{lem}
%Note that when $(q,r) > 2$ the congruence conditions are incompatible.
\begin{proof}
By the Chinese remainder theorem, the pair of simultaneous congruences $n \equiv 0(q)$, $n \equiv -d(r)$ corresponds to the single congruence $n \equiv c([q,r])$, where $c := \frac{q}{(q,r)}r\bar{r} - d\frac{q}{(q,r)}\bar{q} \in \mb{Z}/[q,r]\mb{Z}$; here, $r\bar{r} \equiv 1 (q/(q,r))$ and $\frac{q}{(q,r)}\bar{q} \equiv 1 (r)$. Note that $[q,r]$ is squarefree, and moreover $p|([q,r],c)$ if, and only if, $p|q$, and similarly, $p|([q,r],c+d)$ if, and only if, $p|r$. Applying the previous lemma with these observations implies \eqref{PLAIN}. \\
For \eqref{NONPLAIN}, we note that $[q,r,e] = qre$ whenever $(n,a) = 1$ and $n \equiv 0 (q)$. Consequently,
\begin{align*}
&\sum_{n \leq x \atop n \equiv 0 (q), n+a \equiv 0 (r)} \mu^2(n)\mu^2(n+a)1_{(n,a) = 1} = \sum_{e|a} \mu(e)\sum_{n \leq x \atop n \equiv 0 (qe), n+a \equiv 0(r)} \mu^2(n)\mu^2(n+a) \\
&=E(x)\prod_{p_1^2|| a}\left(1+\frac{1}{p_1^2}\right) \sum_{e|a} \mu(e) \prod_{p_2|qer} \frac{1}{p_2}\left(1-\frac{1}{p_2}\right)\left(1-\frac{2}{p_2^2}\right)^{-1} + O\left(x^{\frac{2}{3}+\e}\tau(qr)\tau(a)\right) \\
&=E(x)\prod_{p_1^2|| a}\left(1+\frac{1}{p_1^2}\right) \prod_{p_2|qr}\frac{1}{p_2}\left(\frac{1-1/p_2}{1-2/p_2^2}\right) \sum_{e|a} \mu(e) \prod_{p_3|e} \frac{1}{p_3}\left(\frac{1-1/p_3}{1-2/p_3^2}\right) + O\left(x^{\frac{2}{3}+\e}\tau(qr)\tau(a)\right)\\
&= E(x;a)\prod_{p|qr}\frac{1}{p}\left(1-\frac{1}{p}\right)\left(1-\frac{2}{p^2}\right)^{-1} +O_a\left(x^{\frac{2}{3}+\e}\tau(qr)\right),
\end{align*}
as claimed.
\end{proof}
We may now prove that the set of squarefree integers is siftable and regular.
\begin{proof}[Proof of Proposition \ref{QUADRATFREI}]
Referring to Definition \ref{SIFTABLE}, i) and ii) hold by Lemma \ref{SIEVESTUFF} where we clearly see that $R(x;q,r) \ll x^{\frac{2}{3} + \e} \tau(qr)f(qr)^{-1}$, with 
\begin{equation} \label{FUPBD}
f(p) := \frac{1}{p}\left(1-\frac{1}{p}\right)\left(1-\frac{2}{p^2}\right)^{-1} \leq \frac{1}{p}. 
\end{equation}
Moreover, note that if $\theta \in (0,1/6-2\e)$ then if $qr \leq x^{\theta}$, 
\begin{equation*}
\sum_{qr \leq x^{\theta}} R(x;q,r) \ll x^{\frac{2}{3}+\e} \sum_{qr \leq x^{\theta}} (qr)^{1+\e} \ll x^{\frac{2}{3} + 2\theta + 2\e} \ll x^{1-2\e},
\end{equation*}
so iii) holds as well with $\theta = 1/7$, say. Note that \eqref{FUPBD} and Mertens' theorem implies that $\sum_{p \leq x} f(p) = \log_2 x + O(1)$, so $S$ is also regular.
\end{proof}
Theorem \ref{CHOW} now follows.
\begin{proof}[Proof of Theorem \ref{CHOW}]
We first observe that $(n,a) = (n+a,a)$. As such, we can write $\mu_y(n)\mu_y(n+a) = \mu_y(n/(n,a))\mu_y((n+a)/(n,a))$, so that
\begin{equation} \label{DIVCOMB}
\sum_{n \leq x} \mu_y(n) \mu_y(n+a) = \sum_{\alpha | a}\sum_{m \leq x/\alpha} \mu_y(m)\mu_y(m+a/\alpha)1_{(m,a/\alpha) = 1}.
\end{equation}
%Given divisors $\alpha_1,\alpha_2|a$ we can partition the set of pairs $\{(n,n+a) : n \leq x\}$ into sets $F_{\alpha_1,\alpha_2} := \{n \leq x : (n,a) = \alpha_1, (n+a,a) = \alpha_2\}$. Note that 
%Then observe that
%\begin{equation*}
%\sum_{n \leq x} \mu(n) \mu(n+a) = \sum_{\alpha_1,\alpha_2| a} \sum_{n \in F_{\alpha_1,\alpha_2}} 
%x
Now, if $\alpha$ is fixed $b := a/\alpha$ and $x' := x/\alpha$ then
\begin{align*}
&\left|\sum_{n \leq x'} \mu_y(n)\mu_y(n+b)1_{(n,b) = 1}\right| = \left|\sum_{n \leq x} \mu^2(n)\mu^2(n+b)1_{(n,b) = 1} e^{\pi i\omega_y(n)}e^{\pi i \omega_y(n+b)}\right| \\
&= \left|\sum_{n \leq x'} \mu^2(n)\mu^2(n+b)1_{(n,b) = 1} e^{\pi i\sqrt{\log_2 y}\tilde{\omega}_y(n)}e^{\pi i\sqrt{\log_2 y}\tilde{\omega}_y(n+b)}\right| \\
&= E(x';b)|\phi_{x',y}(\pi\sqrt{\log_2 y},\pi\sqrt{\log_2 y})|.
\end{align*}
In light of Proposition \ref{QUADRATFREI}, we may apply Proposition \ref{PHIEST}, yielding
\begin{equation*}
x^{-1}\sum_{n \leq x'} \mu_y(n)\mu_y(n+b)1_{(n,b) = 1} \ll_a \frac{1}{\log_2 y} + e^{-\frac{1}{21} \beta \log \beta} + \log^{-\frac{1}{6}} x. 
%= \frac{1}{\log(\beta \log x)} + e^{-\frac{1}{6}\beta \log \beta}.
\end{equation*}
Inserting this into \eqref{DIVCOMB} for each of the finitely many divisors $\alpha$ of $a$ completes the proof of \eqref{ONLYZ}. Equation \eqref{SMALLU} follows immediately from \eqref{B4} and the triangle inequality.
%\begin{equation*}
%\phi_x(u,v) 
%\begin{proof}[Proof of Theorem \ref{EKTHM}]
%By Lemma \ref{ESSEEN}, we have
%\begin{equation*}
%\|F_x-\Phi\|_{L^{\infty}(\mb{R}^2)} 
\end{proof}
\section{A Disjunction Theorem for Characteristic Functions}
In this section, we show that if the distance between the characteristic function of a bivariate Gaussian random vector and that of a given random vector is not small then the $L^{\infty}$ distance between the corresponding distribution functions of these random vectors is ''smaller than expected''. In the final section of this paper we will leverage this fact to gain better control on the difference between $\phi_{x,y}(u,v)$ and $\chi(u,v)$, in the situation that $S$ is the set of squarefree integers. To begin with, we will need the following arithmetic estimate that is relevant to this situation.
\begin{lem}\label{SPIRO}
Let $C > 0$. Let $k \in \mb{N}$ such that $|k-\log_2 x| \leq C \sqrt{\log_2 x}$. Then
\begin{equation*}
\tilde{\pi}_k(x) := \left|\{n \leq x : \mu^2(n(n+1)) = 1, \omega(n) = k\}\right| \gg_C \frac{x}{\sqrt{\log_2 x}}.
\end{equation*}
\end{lem} 
\begin{proof}
Let $Y \geq 2$. Given a modulus $q$ and a residue class $a$ modulo $q$, write $\pi_k(x;a,q)$ to denote the number of integers $n \leq x$ such that $\omega(n) = k$ and $n \equiv a(q)$. Then we have
\begin{align*}
\tilde{\pi}_k(x) &= \sum_{n \leq x \atop \omega(n) = k} \mu^2(n(n+1)) = \sum_{e_1e_2 \leq Y \atop (e_1,e_2) = 1} \mu(e_1)\mu(e_2) \pi_k(x;(e_1e_2)^2,c(e_1,e_2)) \\
&+ O\left(\sum_{e_1,e_2 \leq x^{\frac{1}{2}} \atop e_1e_2 > Y} \pi_k(x;(e_1e_2)^2,c(e_1,e_2))\right) \\
&=: T_1 + T_2,
\end{align*}
where $c(e_1,e_2)$ is the residue class modulo $(e_1^2e_2^2)$ corresponding to the pair of congruences $n \equiv 0 (e_1^2)$ and $n \equiv -1 (e_2^2)$. Put $\rho := \frac{k-1}{\log_2 x}$. Fix $e_1$ and $e_2$ momentarily, and let $q := (e_1e_2)^2$. By Theorem 2 of \cite{SPI}, the asymptotic formula
\begin{equation*}
\pi_k(x;q,c(e_1,e_2)) = \frac{1}{\phi(q)}\frac{x (\log_2 x)^{k-1}}{(k-1)!\log x} \left(\frac{\alpha(\rho;q)}{\Gamma\left(1+\rho\right)}\left(\frac{\phi(q)}{q}\right)^{\rho} + O_{\tau}\left(\frac{\rho}{\log_2 x}(\log_4 x)^2\right)\right)
%\left[\frac{\partial^{k-1}}{\partial z^{k-1}} \frac{1}{\Gamma(z+1)}\alpha (z)\prod_{p|q}\left(1+\frac{z}{p-1}\right)\log^z x\right]_{z = 0} \left(1+O_{\tau}\left(\frac{\log_2(3q)}{\log x}\right)\right)
\end{equation*}
holds uniformly in $q \leq \log^{\tau} x$, for each fixed $\tau > 0$ and $1 \leq k \leq 2\log_2 x$. Here, we put
\begin{equation*}
\alpha(\rho;q) := \prod_{p\nmid q} \left(1+\frac{\rho(1-\rho)}{p(p-1)} + \rho(1-\rho)h_p(\rho)\right),
\end{equation*}
where the numbers $h_p(\rho)$ satisfy $\sum_p h_p(\rho) \ll 1$.
%Here, $\alpha(z)$ is the residue of the product $(s-1)^z\prod_{p}\left(1+\frac{z}{p^s-1}\right)$ at $s = 1$. Put
%\begin{equation*}
%\gamma_k(q) := \left[\frac{\partial^{k-1}}{\partial z^{k-1}} \frac{1}{\Gamma(z+1)}\alpha (z)\prod_{p|q}\left(1+\frac{z}{p-1}\right)\log^z x\right]_{z = 0} =:\left[\frac{\partial^{k-1}}{\partial z^{k-1}} \gamma(z;q)e^{z\log_2 x}\right]_{z = 0}.
%\end{equation*}
%Of course, for $\text{Re}(s) > 1$, $(s-1)^z\prod_p\left(1+\frac{z}{p^s-1}\right) \asymp \prod_p \left(1+\frac{z}{p^s-1}\right)\left(1-\frac{1}{p^s}\right)^{-z}$ is a holomorphic function in a neighbourhood of $s = 1$, and is independent of $x$. Thus, $\alpha(z)$ and its derivatives are bounded in $x$. The same is true regarding the contribution to the derivatives from $1/\Gamma(z+1)$. Finally, as $q < \log^{\tau} x$, we have 
%\begin{equation*}
%\lambda_q(z) := \prod_{p|q}\left(1+\frac{z}{p-1}\right)= f(z)\exp\left(z\sum_{p | q} \frac{1}{p-1}\right),
%\end{equation*}
%where $f(z)$ is holomorphic in $z$ and uniformly bounded as $q \ra \infty$. Thus, as $\lambda_q(0) = f(0) = 1$, for each $t \leq k$,
%\begin{equation*}
%\left[\left(\frac{\partial}{\partial z}\right)^{t-1} \lambda_q(z)\right]_0 = \lambda_q(0)\sum_{l \leq t} \binom{t}{l} f^{(t-l)}(0)\left(\sum_{p|q} \frac{1}{p-1}\right)^l \gg \left(\sum_{p|q} \frac{1}{p-1}\right)^t.
%\end{equation*}
%As such, we have
%\begin{equation*}
%\gamma_k = \sum_{0 \leq l \leq k-1} \binom{k-1}{l} \gamma^{(l)}(0;q) \log_2^{k-l} x = 
Choose $Y = \log^6 x$ and put $\rho = 1+O\left(\frac{1}{\sqrt{\log_2 x}}\right)$. Since 
\begin{equation*}
\alpha(\rho;q) = \exp\left(O\left(\rho(1-\rho)\sum_{p\nmid q}\frac{1}{p(p-1)}\right)\right) = 1 + O\left(|\rho-1|\right) = 1+O\left(\frac{1}{\sqrt{\log_2 x}}\right),
\end{equation*}
and also
\begin{equation*}
(\phi(q)/q)^{\rho} = \frac{\phi(q)}{q}\left(1+O\left( \frac{\log_2 q}{\sqrt{\log_2 x}}\right)\right) = \frac{\phi(q)}{q}\left(1 + O\left(\frac{\log_3 x}{\sqrt{\log_2 x}}\right)\right),
\end{equation*}
we have
\begin{align*}
T_1 &= \frac{x(\log_2 x)^{k-1}}{(k-1)! \log x} \left(\sum_{(e_1,e_2) = 1} \frac{\mu(e_1)\mu(e_2)}{e_1^2e_2^2} + O\left(\frac{\log_3 x}{\sqrt{\log_2 x}}\right)\right) \\
&= \left(1+O\left(\frac{\log_3 x}{\sqrt{\log_2 x}}\right)\right) \frac{E(x)(\log_2 x)^{k-1}}{(k-1)! \log x}.
\end{align*}
By Stirling's formula, we have
\begin{align*}
\frac{(\log_2 x)^{k-1}}{(k-1)! \log x} &= \left(1+O\left(\frac{1}{\log_2 x}\right)\right) (2\pi)^{-\frac{1}{2}}\left(\frac{e\log_2 x}{k-1}\right)^{k-1} e^{-\log_2 	x} \frac{1}{\sqrt{\log_2 x}} \\
&= \left(1+O\left(\frac{1}{\log_2 x}\right)\right)(2\pi)^{-\frac{1}{2}}e^{k-1-\log_2 x} \cdot \left(1+\frac{k-1-\log_2 x}{\log_2 x}\right)^{-(k-1)} \frac{1}{\sqrt{\log_2 x}} \\
&\gg_C \frac{1}{\sqrt{\log_2 x}},
\end{align*}
so that $T_1 \gg_C x(\log_2 x)^{-\frac{1}{2}}$. It now suffices to show that $T_2 = o\left(\frac{x}{\sqrt{\log_2 x}}\right)$ to complete the proof. \\
For the sum in $T_2$, we simply bound the terms $\pi_k(x;(e_1e_2)^2,c(e_1,e_2))$ by the number of integers $n \leq x$ satisfying $\mu^2(n(n+1) = 1$ with $n \equiv c (e_1^2e_2^2)$ and apply the bound in the proof of Lemma \ref{HB} (with $q = 1$ there), getting
\begin{equation*}
T_2 \ll xY^{-\frac{1}{2}} \log^2 x \ll x \log^{2-\tau/2} x \ll \frac{x}{\log x}.
\end{equation*}
%Hence, we have
%\begin{equation*}
%\tilde{\pi}_k(x) = C\frac{x}{\log x} \left(\frac{\gamma_k}{(k-1)!}\left(1+O\left(\frac{\log_3 x}{\log x}\right)\right) + O\left(\frac{1}{\log x}\right)\right).
%\end{equation*}
%Specializing to $k=(1+O(1/\sqrt{\log_2 x}))\log_2 x$, it suffices to show that $\gamma_k \gg \log_2^{k-1} x$. 
The claim follows.
\end{proof}

We surmise from the Lemma \ref{SPIRO} the minimum order of growth of the $L^{\infty}$ distance between $F_{x,y}$ and $\Phi_{(2)}$, where $F_{x,y}$ is the distribution function whose characteristic function is $\phi_{x,y}$.
\begin{prop}\label{CONTRA}
Let $C > 0$ be fixed. Let $y := x^{\frac{1}{\beta}}$ with $\beta \leq e^{B(x)}$. Let $\mf{L}_y := \sqrt{\log_2 y}$, 
\begin{equation*}
F_{x,y}(t,t') := E(x)^{-1}\left|\left\{n \leq x: \mu^2(n(n+1)) = 1, \tilde{\omega}_y(n) \leq t, \tilde{\omega}_y(n+1) \leq t'\right\}\right|.
\end{equation*}
Then the estimate
\begin{equation}
\|F_{x,y} - \Phi_{(2)}\|_{L^{\infty}(\mb{R}^2)} = o\left(\mf{L}^{-1}\right) \label{AST}
\end{equation}
cannot hold.
\end{prop}
\begin{proof}
Assume that \eqref{AST} holds. Let $t = \log_2 y$. By the Tur\'{a}n-Kubilius inequality (see Section III.2 in \cite{Ten2}), for any $t' \in \mb{R}$ we have
\begin{equation*}
\left|F_{x,y}(t',t)-F_{x,y}(t',\infty)\right| \ll x^{-1} |\{n \leq x : |\omega_y(n+1)-\log_2 y| > \log_2 y\}| \ll \frac{1}{\log_2 y}.
\end{equation*}
Put now $t' = 0$. Since $\Phi(t) \ll 1/\log x$,
\begin{align*}
&F_{x,y}(t'+(2\mf{L}_y)^{-1},\infty)-F_{x,y}(t'-(2\mf{L}_y)^{-1},\infty) \\
&= F_{x,y}(t'+(2\mf{L}_y)^{-1},t)-F_{x,y}(t'-(2\mf{L}_y)^{-1},t) + O\left(\mf{L}_y^{-2}\right) \\
&= \Phi(t)\left(\Phi(t'+(2\mf{L}_y)^{-1})-\Phi(t'-(2\mf{L}_y)^{-1})\right) + O\left(\mf{L}^{-2} + \|F_{x,y}-G\|_{L^{\infty}(\mb{R}^2)}\right) = o\left(\mf{L}^{-1}\right).
\end{align*}
On the other hand, by Lemma \ref{SPIRO},
\begin{align*}
&F_{x,y}(t'+(2\mf{L}_y)^{-1},\infty)-F_{x,y}(t'-(2\mc{L}_y)^{-1},\infty) \\
&=  E(x)^{-1}|\{n \leq x: \mu^2(n(n+1)) = 1, \left|\omega_y(n) - \log_2 y\right| < 1 \}|\\
&\geq E(x)^{-1}|\{n \leq x: \mu^2(n(n+1)) = 1, \left|\omega(n) - \log_2 x\right| < 1 \}| \gg \mf{L}^{-1}.
\end{align*}
%upon recalling that $\log B(x) = o(\sqrt{\log_2 x})$. 
This is an obvious contradiction.
\end{proof}
The last two lemmata imply the following disjunction result, which allows one to extract additional savings on the distance between $\phi_{x,y}$ and the characteristic function of the bivariate Gaussian distribution.
\begin{prop} \label{MAINPROP}
Let $Z,T \geq 1$, $y = x^{\frac{1}{\beta}}$ with $\beta \leq e^{B(x)}$ and $u,v \in \mb{R}$. Let $F_{x,y}$, $\Phi_{(2)}$, $\phi_{x,y}$, $\chi$ and $\lambda(y)$ be as above.
%, with respective characteristic functions $\phi_{x,y}$ and $\chi$. Also, put $\lambda(y) := \sum_{p \leq y} \frac{1}{p}$. 
Then exactly one of the following holds: \\
i) we have
\begin{equation*}
|\phi_{x,y}(u,v) - \chi(u,v)| \ll Z\left(\frac{1}{T} + \int_{[T^{-3},T]^2} \frac{\left|\phi_{x,y}(t,t') - \chi(t,t')\right|}{|tt'|}dtdt' \right);
\end{equation*}
ii) there is a bivariate Poisson random vector $X = (X_1,X_2) \sim \text{Poi}(\lambda(y))^2$ defined on $[1,x]$ with $X_1,X_2$ independent, such that if $\tilde{X}_j := (X_j-\lambda(y))/\sqrt{\lambda(y)}$ and $\tilde{X} = (\tilde{X}_1,\tilde{X}_2)$ then
\begin{equation*}
\|F_{x,y}-\Phi_{(2)}\|_{L^{\infty}(\mb{R}^2)} \ll Z^{-1}\left(d_{TV}(\mc{L}(\tilde{\omega}_y,\tilde{\omega}_y \circ T),\mc{L}(\tilde{X}))+e^{-\frac{1}{2}w^2}w^4(\log_2 y)^{-1}\right).
\end{equation*}
%where $\pi_x$ is the characteristic function of $\tilde{X}_1$, and $\chi$ is the characteristic function of a normal random variable.
\end{prop}
\begin{proof}
Suppose first that $|\phi_{x,y}(u,v)-\chi(u,v)| \leq Z\|F_{x,y}-\Phi_{(2)}\|_{\infty}$. Then i) follows upon applying Lemma \ref{ESSEEN} together with Corollary \ref{APPLY}. \\
Suppose now that $\|F_{x,y}-G\|_{\infty} < Z^{-1}|\phi_{x,y}(u,v)-\chi(u,v)|$. Write $\lambda$ in place of $\lambda(z)$ for convenience, and let $X = (X_1,X_2)$ be as in the statement. Thus,
\begin{align*}
\|F_{x,y}-\Phi_{(2)}\|_{\infty} &< Z^{-1}\left(|\phi_{x,y}(u,v)-\pi_y(u)\pi_y(v)| + |\pi_y(u)||\pi_y(v)-\chi(v)| + |\chi(v)||\pi_y(u)-\chi(u)|\right) \\
&\ll Z^{-1}\left(|\phi_{x,y}(u,v)-\pi_y(u)\pi_y(v)| + \max_{w \in \{u,v\}} |\pi_y(w)-\chi(w)|\right).
\end{align*}
By Lemma \ref{TRIV} , $|\phi_{x,y}(u,v)-\pi_{y}(u)\pi_y(v)| \ll d_{TV}(\mc{L}(\tilde{\omega}_y,\tilde{\omega}_y \circ T), \mc{L}(\tilde{X}))$. Furthermore, the proof of Proposition \ref{PHIEST} shows that
%observe that 
%\begin{equation*}
%\pi_z(w) = \sum_{k \geq 0} e^{iw\frac{k-\lambda}{\sqrt{\lambda}}} \mb{P}(X_1 = k) = e^{-\lambda - iw\sqrt{\lambda}}\sum_{k \geq 0} e^{i\frac{w}{\sqrt{\lambda}}k} \frac{\lambda^k}{k!} = e^{-iw\sqrt{\lambda}} e^{\lambda\left(e^{iw/\sqrt{\lambda}}-1\right)}.
%\end{equation*}
%$|\pi_y(w)| = \exp\left(-\lambda(1-\cos(w/\sqrt{\lambda}))\right)$. Clearly, $\chi(w) = e^{-\frac{1}{2}w^2}$, so that
\begin{equation*}
|\pi_y(w)-\chi(w)| \leq e^{-\frac{1}{2}w^2}w^4(\log_2 y)^{-1}.
\end{equation*}
This implies the claim.
\end{proof}
\section{On a Problem of Erd\H{o}s and Mirsky}
%We may briefly explain (our version of) the Bateman-Spiro heuristic here using basic ideas from Section 6. 
We begin this final section by explaining a version of the Bateman-Spiro heuristic underlying the statement of Conjecture \ref{CONJEPS}. Recall that $$S(x) := \{n \leq x: \tau(n) = \tau(n+1)\}.$$ In light of the result in \cite{EPS} it suffices to show that $S(x) \gg x(\log_2 x)^{-\frac{1}{2}}$. Let 
\begin{equation*}
S^{\ast}(x) := |\{n \leq x : \mu^2(n(n+1)) = 1, \tau(n) = \tau(n+1)\}|.
\end{equation*}
Clearly, $S(x) \geq S^{\ast}(x)$, and we seek to bound $S^{\ast}(x)$ from below. 
By Lemmata \ref{MELLINLEM} and \ref{GX} below (neither of which are deep) we have
\begin{equation} \label{STATE}
S^{\ast}(x) \gg \frac{E(x)}{\log_2 x}\left|\int_0^{2\pi\sqrt{\log_2 x}}\int_0^{2\pi\sqrt{\log_2 x}} \phi_x(u-u',u'-u)dudu'\right|.
\end{equation}
Suppose that we had some control over the error term $|\phi_x(v,-v)-\chi(v,-v)|$ in such a way that this term is negligible for each $v \in [-2\pi \sqrt{\log_2 x},2\pi \sqrt{\log_2 x}]$.  Then we could replace $\phi_x(u-u',u'-u)$ in \eqref{STATE} by $e^{-(u-u')^2}$. It is then easy to prove (see Lemma \ref{HEUR}) that the corresponding integral is $\gg \sqrt{\log_2 x}$. When used in \eqref{STATE}, this gives the conjectured lower bound $\frac{x}{\sqrt{\log_2 x}}$ for $S^{\ast}(x)$ \emph{on heuristic grounds}. \\
To furnish the lower bound rigorously, we would need to show that we can suitably bound $|\phi_x(u-u',u'-u)-\chi(u-u',u'-u)|$, at least in an average sense, in the box $[0,2\pi \sqrt{\log_2 x}]^2$. As we cannot do this, we settle for the following analogue of Conjecture \ref{CONJEPS}. \\
Let $2 \leq y \leq x$, such that if $\beta := (\log x)/(\log y)$ then $A(x) \leq \beta \leq e^{B(x)}$. We note, in particular, that $\log_2 y = (1+o(1))\log_2 x$. Write $\mf{L}_y:= \sqrt{\log_2 y}$ and $\mf{L} : =\mf{L}_x$ as before. For $|j| \leq \left(1-\e\right)\sqrt{\log_2 x\log_3 x}$ let 
\begin{align*} S_j(x;y) &:= |\{n \leq x: \tau_y(n) = 2^j\tau_y(n+1)\}| \\ S^{\ast}_j(x;y) &:= |\{n \leq x : \mu^2(n(n+1)) = 1, \tau_y(n) = 2^j\tau_y(n+1)\}|, 
\end{align*}
where $\tau_y$ is the $y$-smooth divisor function, i.e., $\tau_y(n) := \prod_{p^k || n \atop p \leq y} (k+1)$.
%Erd\H{o}s and Mirsky \cite{EM} famously conjectured that $\lim_{x \ra \infty} S(x) = \infty$, i.e., that there are infinitely many integers $n$ such that $\tau(n) = \tau(n+1)$. Based on ideas of Spiro, Heath-Brown proved this conjecture in the affirmative, proving in fact that $S(x) \gg x\log^{-7} x$. This was improved by Hildebrand subsequently to $S(x) \gg x\log_2^{-3} x$, based on an upper bound argument due to Erd\H{o}s, Pomerance and Sark\H{o}zy. The latter proved that $S(x) \ll x\log_2^{-\frac{1}{2}} x$, and conjectured that this upper bound is the correct order of magnitude. They cite a heuristic argument due to Bateman and Spiro as further motivation for this conjecture (though the heuristic of Bateman-Spiro is not explicitly described in their paper). \\
%We will shortly describe this heuristic in detail, and will in fact prove the validity of the heuristic, and thus the conjecture of Erd\H{o}s-Pomerance-Sark\H{o}zy. Let 
 We will estimate $S_j^{\ast}$ from below to give a lower bound for $S_j$.\\
\begin{thm}\label{EM}
Let $2 \leq y \leq x$ such that if $y = x^{\frac{1}{\beta}}$ then $A(x) \leq \beta \leq e^{B(x)}$. 
%Let $\tau_y(n)$ be the number of divisors $d |n$ such that if $p|d$ then $p \leq y$. 
Then for each $|j| \leq (1-\e)\sqrt{\log_2 x \log_3 x}$, we have
\begin{equation*}
|\{n \leq x : \tau_y(n) = 2^j\tau_y(n+1)\}| \gg  e^{-\frac{j^2}{4\mc{L}_y^2}}\frac{x}{\sqrt{\log_2 x}}.
\end{equation*}
\end{thm}
\begin{rem}
Note that the set of $n \leq x$ such that $\tau(n) = \tau(n+1)$ distributes among sets in which $\tau_y(n) = 2^j \tau_y(n+1)$, for some $|j| \leq \frac{\log x}{\log y}$. However, most integers have about $\log\left(\frac{\log x}{\log y}\right)$ such factors, and we thus expect that $n$ and $n+a$ should have nearly the same number of prime factors of size larger than $y$. Thus, the restriction on $j$ here is not expected to be problematic when $y$ is close to the upper limit of its range. If we had a better understanding of which $j$ occur most often such that $\tau(n) = \tau(n+1)$ \emph{and} $\tau_y(n) = 2^j\tau_y(n+1)$ then we would have a chance at actually proving Conjecture \ref{CONJEPS}.
%\\ 
%Our argument also suffices to deduce the upper bound $
\end{rem}
%It is clear that $S^{\ast}(x;y) \leq S(x;y)$, and we shall estimate $S^{\ast}(x;y)$ from below in order to get a lower bound for $S(x;y)$. \\
The relevance of the results of the previous sections to this problem is brought to light via the following lemma. 
\begin{lem} \label{MELLINLEM}
Let $x \geq 3$ and $\alpha > 0$. Then
\begin{equation*}
S^{\ast}_j(x;y) = \frac{E(x)}{(\pi\mf{L}_y\log 2 )^2} \int_{[0,2\pi \mf{L}_y]^2}g_{y,\alpha}(u)g_{y,\alpha}(u')\phi_{x,y}(u-u',u'-u)e^{i(j/\mf{L}_y)(u'-u)}d\mbf{u},
%\sum_{n \leq x} \mu^2(n(n+1)) e^{i(t-t')\tilde{\omega}_x(n)}e^{i(t'-t)\tilde{\omega}_x(n+1)},
%\frac{1}{(\pi\log 2)^2} \int_{\mb{R}^2} \frac{dt dt'}{(\sg + it/\log 2)(\sg+it'/\log 2)} \sum_{n \leq x} \mu^2(n(n+1)) e^{i(t-t')\omega(n)}e^{i(t'-t) \omega(n+1)}.
\end{equation*}
where we set
\begin{align*}
g_{y,\alpha}(u) &:= \lim_{M \ra \infty} \sum_{|k| \leq M} \frac{1}{\alpha + i\left(2\pi k + u\right)/\log 2} \\
&= \frac{1}{\alpha + iu/\mf{L}_y} + \sum_{k \geq 1} \frac{\alpha + iu/(\mf{L}_y\log 2)}{(\alpha + iu/(\mf{L}_y\log 2))^2 + 4\pi^2k^2/\log^2 2}.
\end{align*}
\end{lem}
\begin{proof}
Let $b > 0$. Applying the Mellin identity 
\begin{equation}
\frac{1}{2\pi i}\int_{(\alpha)} b^s\frac{ds}{s} = \lim_{T \ra \infty} \left(\frac{1}{2\pi i} \int_{\alpha-iT}^{\alpha+iT} b^s \frac{ds}{s}\right)= \begin{cases} 1 &\text{ if $b > 1$} \\ \frac{1}{2} &\text{ if $b=1$} \\ 0 &\text{ if $b < 1$}\end{cases} \label{MELLIN}
\end{equation}
for $b = \tau(n)/(2^j\tau(n+1))$ and its reciprocal yields
\begin{equation*}
(\pi i)^{-2}\lim_{T,T' \ra \infty} \left(\int_{\alpha-iT}^{\alpha+iT}\int_{\alpha-iT'}^{\alpha+iT'} \frac{dsds'}{ss'} \left(\frac{\tau(n)}{2^j\tau(n+1)}\right)^s \left(\frac{2^j\tau(n+1)}{\tau(n)}\right)^{s'}\right) = \begin{cases} 1 &\text{ if $\frac{\tau(n)}{\tau(n+1)} = 2^j$} \\ 0 &\text{ otherwise}. \end{cases}
\end{equation*}
Summing this identity over squarefree integers $n \leq x$, we have
\begin{align*}
&S_j^{\ast}(x;y) = \frac{1}{\pi^2} \lim_{T,T' \ra \infty} \left(\int_{-T}^{T}\int_{-T'}^{T'} \frac{dtdt'}{(\alpha + it)(\alpha+it')} \sum_{n \leq x \atop \mu^2(n(n+1)) = 1}  \tau_y(n)^{i(t-t')} \tau_y(n+1)^{i(t'-t)}2^{ij(t'-t)}\right) \\
&= \frac{1}{\pi^2} \lim_{T,T' \ra \infty} \left(\int_{-T}^{T}\int_{-T'}^{T'} \frac{dtdt'}{(\alpha + it)(\alpha+it')} \sum_{n \leq x \atop \mu^2(n(n+1)) = 1} e^{i\omega_y(n)(t-t')\log 2 } e^{i\omega_y(n+1)(t'-t)\log 2 } e^{i(t'-t) j\log 2}\right) \\
&= \frac{1}{(\pi \log 2)^2} \lim_{T,T' \ra \infty} \left(\int_{-T}^{T}\int_{-T'}^{T'} \frac{dtdt'}{(\alpha + \frac{it}{\log 2})(\alpha+\frac{it'}{\log 2})} \sum_{n \leq x \atop \mu^2(n(n+1)) = 1} e^{i(t-t') \omega_y(n)} e^{i(t'-t)\omega_y(n+1)}e^{ij(t'-t)}\right).
\end{align*}
We now observe that the sum over $n$ is invariant as a function of $t$ and $t'$ under the translations $(t,t') \mapsto (t-2\pi k,t-2\pi l)$, for $k, l \in \mb{Z}$. Since the double limit defining $S_j^{\ast}(x,y)$ exists, we may rewrite it as
\begin{align*}
S_j^{\ast}(x;y) &= \frac{1}{(\pi \log 2)^2}\lim_{M,N \ra \infty \atop M,N \in \mb{Z}} (\int_{-2\pi M}^{2\pi M} \int_{-2\pi N}^{2\pi N}\frac{dtdt'}{(\alpha + it/\log 2)(\alpha+it'/\log 2)} \\
&\cdot \sum_{n \leq x \atop \mu^2(n(n+1)) = 1}  e^{i(t-t') \omega_y(n)} e^{i(t'-t)\omega_y(n+1)}e^{ij(t'-t)}) \\
&=\frac{1}{(\pi \log 2)^2}\int_{[0,2\pi]^2} \left(\sum_{n \leq x} \mu^2(n(n+1)) e^{i(t-t') \omega_y(n)} e^{i(t'-t)\omega_y(n+1)}e^{ij(t'-t)}dtdt'\right) \\
&\cdot \lim_{M,N \ra \infty \atop M,N \in \mb{Z}} \left(\sum_{|k| \leq M} \frac{1}{\alpha + i(t+2\pi k)/\log 2}\right)\left(\sum_{|l| \leq N} \frac{1}{\alpha+i(t'+2\pi l)/\log 2}\right).
\end{align*}
We may simplify the sum over $k$ via
\begin{align*}
\sum_{|k| \leq M} \frac{1}{\alpha + i(t+2\pi k)/\log 2} &= \frac{1}{\alpha + it/\log 2} + \sum_{1 \leq k \leq M} \left(\frac{1}{\alpha + i(t+2\pi k)/\log 2} + \frac{1}{\alpha + i(t-2\pi k)/\log 2}\right) \\
&= \frac{1}{\alpha + it/\log 2} + 2\sum_{1 \leq k \leq M} \frac{\alpha + it/\log 2}{(\alpha + it/\log 2)^2+4\pi^2 k^2/\log^2 2}.
\end{align*}
A similar expression exists for the sum over $l$ (as a function of $t'$ instead of $t$). Since each of these sums converges absolutely and uniformly in $[0,2\pi]$, we may take $M,N \ra \infty$ to give
\begin{align*}
S_j^{\ast}(x;y) &= \frac{1}{(\pi \log 2)^2}\int_{[0,2\pi]^2} \left(\sum_{n \leq x \atop \mu^2(n(n+1))=1} e^{i(t-t')\mf{L}_y\tilde{\omega}_y(n)} e^{i(t'-t)\mf{L}_y\tilde{\omega}_y(n+a)}\right)e^{ij(t'-t)}g_{y,\alpha}(t\mf{L}_y)g_{y,\alpha}(t'\mf{L}_y)d\mbf{t}.
%\left(\frac{1}{\alpha + it/\log 2} + 2\sum_{k \geq 1} \frac{\alpha + it/\log 2}{(\alpha + it/\log 2)^2+4\pi^2 k^2/\log^2 2}\right)\left(\frac{1}{\alpha + it'/\log 2} + 2\sum_{l \geq 1} \frac{\alpha + it'/\log 2}{(\alpha + it'/\log 2)^2+4\pi^2 l^2/\log^2 2}\right).
\end{align*}
Making the change of variables $u := t\mf{L}_y$ amd $u' := t'\mf{L}_y$ completes the proof of the lemma.
\end{proof}
We will choose $\alpha$ suitably so that $g_{y,\alpha}$ is well-behaved. The following lemma gives us a hint in this direction.
\begin{lem}\label{GX}
Uniformly in $[0,2\pi \mf{L}_y]$, we have $g_{y,\alpha}(u) = \frac{1}{4} + O\left(\frac{1}{\alpha}\right)$, as $\alpha \ra \infty$.
\end{lem}
\begin{proof}
Note that $(\alpha + iu/\mf{L}_y)^{-1} = \frac{1}{\alpha}\left(1+O\left(\frac{1}{\alpha}\right)\right)$ uniformly, since $|u/\mf{L}_y|$ is uniformly bounded.  For the series, we have
\begin{align*}
&\sum_{k \geq 1}\frac{\alpha+iu/(\mf{L}_y\log 2)}{\alpha^2 -u^2/(\mf{L}_y\log 2)^2 + 4\pi^2k^2 +2i\alpha u/(\mf{L}_y\log 2)} \\
&= \alpha\left(1+O\left(\frac{1}{\alpha}\right)\right)\sum_{k \geq 1} \frac{1}{4\pi^2k^2 + (\alpha^2-u^2/(\mf{L}_y\log 2)^2)} \\
%\left(1+O\left(\frac{1}{\alpha}\right)\right) \\
&= \alpha \left(1+O\left(\frac{1}{\alpha}\right)\right)\sum_{k \geq 1} \frac{1}{4\pi^2k^2 + \alpha^2}.
\end{align*}
%We split the sum over $k$ into the ranges $k \leq \alpha/2\pi$ and its complement. In the first range,
Approximating the latter sum by an integral, we have
\begin{align*}
&\sum_{k \geq 1} \frac{1}{4\pi^2k^2 + \alpha^2} = \frac{1}{2\pi\alpha} \left(\frac{2\pi}{\alpha}\sum_{k \geq 1} \frac{1}{1 + \left(\frac{2\pi k}{\alpha}\right)^2}\right) \\
&= \frac{1}{2\pi \alpha} \int_0^{\infty} \frac{dt}{1+t^2} + \frac{1}{2\pi \alpha}\left(\frac{2\pi}{\alpha} \sum_{k \geq 1} \int_k^{k+1} \left(\frac{1}{1 + \left(\frac{2\pi k}{\alpha}\right)^2} - \frac{1}{1+\left(\frac{2\pi t}{\alpha}\right)^2}\right)dt\right) \\
&= \frac{1}{4\alpha} + O\left(\sum_{k \geq 1} \frac{k}{(\alpha^2+4\pi^2k^2)^2}\right).
\end{align*}
We estimate the error term by splitting the sum into the ranges $k \leq \alpha/2\pi$ and its complement. In the first range,
\begin{equation*}
\sum_{k \leq \alpha/2\pi} \frac{k}{(\alpha^2+ 4\pi^2k^2)^2} \ll \frac{1}{\alpha^4} \sum_{k \leq \alpha/2\pi} k \ll \frac{1}{\alpha^2}.
\end{equation*}
In the second range,
\begin{equation*}
\sum_{k > \alpha/2\pi} \frac{k}{(\alpha^2 + 4\pi^2k^2)^2}  \ll \sum_{k > \alpha/2\pi} \frac{1}{k^3} \ll \frac{1}{\alpha^2}.
\end{equation*}
%
%&\geq \frac{1}{\alpha^2} \sum_{k \leq \alpha/2\pi} 1 \gg \frac{1}{2\pi \alpha} \\
%
%
%
%\sum_{k > \alpha/2\pi} \frac{1}{4\pi^2k^2 + \alpha^2} &\geq \frac{1}{8\pi^2}\sum_{k > \alpha/2\pi} \frac{1}{k^2} \gg \frac{1}{\alpha}.
%\end{align*}
It follows from this that $g_{y,\alpha}(u) = \frac{1}{4} + O\left(\frac{1}{\alpha}\right)$ uniformly in the interval, and the proof is complete.
%, and $C(\alpha) \gg 1$ for $\alpha$ sufficiently large (independently of $x$).
\end{proof}
As a consequence, we have the following.
\begin{lem}\label{HEUR}
There is an $\alpha_0 > 0$ such that for $\alpha \geq \alpha_0$,
\begin{equation}
\int_{[0,2\pi \mf{L}_y]^2} g_{y,\alpha}(u)g_{y,\alpha}(u') e^{i(j/\mf{L}_y)(u'-u)}e^{-(u-u')^2}d\mbf{u} = \left(\frac{\pi^{\frac{3}{2}}}{8} + O\left(\frac{1}{\alpha}\right)\right)e^{-\frac{j^2}{4\mf{L}_y^2}}\mf{L}_y+ O\left(e^{\frac{j^2}{4\mf{L}_y^2}}\right). \label{CONJMAIN}
\end{equation}
\end{lem}
\begin{proof}
By Lemma \ref{GX} we have $g_{y,\alpha}(u)g_{y,\alpha}(u') = \frac{1}{16}+O\left(\frac{1}{\alpha}\right)$ for sufficiently large $\alpha$. Thus, it suffices to show that
\begin{equation*}
\int_{[0,2\pi \mf{L}_y]^2}e^{-(u-u')^2 - i(j/\mf{L}_y)(u-u')}d\mbf{u} = 2\pi^{\frac{3}{2}} e^{-\frac{j^2}{4\mf{L}_y^2}}\mf{L}_y+O\left(e^{\frac{j^2}{4\mf{L}_y^2}}\right).
\end{equation*}
Making a change of variables, the integral on the right may be written as
\begin{equation*}
\int_{[0,2\pi \mf{L}_y]^2} e^{-(u-u')^2-i(j/\mf{L}_y)(u-u')}d\mbf{u}  = e^{-\frac{j^2}{4\mf{L}_y^2}}\int_0^{2\pi \mf{L}_y} du' \left(\int_{-u'}^{u'} + \int_{u'}^{2\pi \mf{L}_y-u'}\right)e^{-\left(v-i\frac{j}{2\mf{L}_y}\right)^2}dv.
\end{equation*}
We note that by standard contour integration,
\begin{equation*}
\int_{-u'}^{u'} e^{-\left(v-i\frac{j}{2\mf{L}_y}\right)^2}dv = \sqrt{\pi} - \int_{|v| > u'} e^{-\left(v-i\frac{j}{2\mf{L}_y}\right)^2}dv = \sqrt{\pi}\left(1 +O\left(e^{-\frac{1}{2}(u')^2}e^{\frac{j^2}{4\mf{L}_y^2}}\right)\right),
\end{equation*}
so that
\begin{equation*}
\int_0^{2\pi \mf{L}_y} du' \int_{-u'}^{u'} e^{-v^2}dv  = 2\pi^{\frac{3}{2}} \mf{L}_y +O\left(e^{j^2/4\mf{L}_y^2}\int_0^{\infty}e^{-\frac{1}{2}(u')^2} du'\right) = 2\pi^{\frac{3}{2}} \mf{L}_y + O\left(e^{\frac{j^2}{4\mf{L}_y^2}}\right).
\end{equation*}
In a similar vein, we have
\begin{equation*}
\int_{u'}^{2\pi \mf{L}_y-u'} e^{-v^2}dv \leq \int_{|v| > u'} e^{-v^2}dv \ll e^{-\frac{1}{2}(u')^2},
\end{equation*}
so that 
\begin{equation*}
\int_0^{2\pi \mf{L}_y} du' \int_{u'}^{2\pi \mf{L}_y-u'} e^{-v^2}dv \ll \int_0^{\infty} e^{-\frac{1}{2}(u')^2}du' \ll 1.
\end{equation*}
This implies the claim.
\end{proof}
%We may now describe the heuristic as follows. The appearance of $\phi_x$ in the identity for $S^{\ast}(x)$ and our previous considerations suggest that $\phi_x$ should be, in an average sense, approximable by $\chi$, the characteristic function of an uncorrelated, bivariate Gaussian random vector. The integral resulting from this approximation bears great resemblance to the integral in the following lemma, which we can estimate quite easily.
%We will show below that the perturbations by the functions $g_{x,\sg}(u)g_{x,\sg}(u')$ are negligible for sufficiently large $\sg$. Thus, the expected main term is described by the previous lemma. 
%In order to turn the heuristic mentioned in the introduction into a proof of Theorem \ref{EM}, we need to estimate the error implicit in replacing $\phi_x(u-u',u'-u)$ by $e^{-(u-u')^2}$. The interest of the methods of the previous sections is in their utility in this direction.
\begin{proof}[Proof of Theorem \ref{EM}] 
As before, we split $$\phi_x(u-u',u'-u) = \chi(u-u',u'-u) + (\phi_x(u-u',u'-u)-\chi(u-u',u'-u)).$$ Thus,
\begin{equation*}
S_j^{\ast}(x;y)/E(x) = \mc{M} + \mc{E},
\end{equation*}
where we put
\begin{align*}
\mc{M} &:= \frac{1}{(\pi \mf{L}_y\log 2)^2}\int_{[0,2\pi\mf{L}_y]^2} g_x(u)g_x(u')e^{-(u-u')^2-i(j/\mf{L}_y)(u-u')}d\mbf{u}, \\
%e^{-4\pi^2 \mf{L}^2\left\|\frac{u-u'}{2\pi \mf{L}}\right\|^2}, \\
\mc{E} &\ll \frac{1}{\mf{L}_y^2} \int_0^{2\pi \mf{L}_y} \int_0^{2\pi \mf{L}_y} \left|E(x)^{-1}\sum_{n \leq x \atop \mu^2(n(n+1)) = 1} e^{i(u-u') \tilde{\omega}_y(n)} e^{i(u'-u)\tilde{\omega}_y(n+1)} - e^{-(u-u')^2}\right|d\mbf{u}.
%- e^{-4\pi^2 \mf{L}^2\left\|\frac{u-u'}{2\pi \mf{L}}\right\|^2}\right|.
\end{align*}
%We will show that $\mc{M} \gg x\log_2^{-\frac{1}{2}} x$. 
By Lemma \ref{HEUR}, $\mc{M} \asymp e^{-\frac{j^2}{4\mf{L}_y^2}}\mf{L}^{-1}$ when $|j| \leq \left(1-\e\right)\sqrt{\log_2 y\log_3 x}$. Our objective now is to show that $\mc{E} = o\left(e^{-\frac{j^2}{4\mf{L}_y^2}}\mf{L}^{-1}\right)$. Making a change of variables as above and using symmetry,
\begin{align*}
\mc{E} &\ll \mf{L}^{-2} \int_0^{\pi \mf{L}_y} du' \int_{-u'}^{2\pi\mf{L}_y- u'} \left|\phi_{x,y}(v,-v) - e^{-v^2}\right|dv \\
&\leq \mf{L}^{-1} \max_{u' \in [0,\pi \mf{L}_y]} \int_{-u'}^{2\pi\mf{L}_y- u'} \left|\phi_{x,y}(v,-v) - e^{-v^2}\right|dv.
\end{align*}
%that by Lebesgue invariance and the $2\pi$-periodicity of the integrand,
%\begin{align}
%\mc{E} &\ll \frac{1}{\mf{L}^2} \int_0^{2\pi \mf{L}}du' \int_{-u'}^{2\pi \mf{L}-u'} dv \left|E(x)^{-1}\sum_{n \leq x} \mu^2(n(n+1)) e^{iv \tilde{\omega}_x(n)} e^{-iv\tilde{\omega}_x(n+1)} -e^{-4\pi^2 \mc{L}^2 \left\|\frac{v}{2\pi \mc{L}}\right\|^2}\right| \nonumber\\
%&\ll \frac{1}{\mc{L}} \int_0^{2\pi \mc{L}} dv\left|E(x)^{-1}\sum_{n \leq x} \mu^2(n(n+1)) e\left(\frac{v}{2\pi \mc{L}}(\omega(n)-\omega(n+1))\right) -e^{-4\pi^2 \mc{L}^2 \left\|\frac{v}{2\pi \mc{L}}\right\|^2}\right|. \label{ERRORS}
%\end{align}
To estimate $\mc{E}$ we use Propositions \ref{CONTRA} and \ref{MAINPROP}. Indeed, combining the proof of the latter proposition with Lemma \ref{TVPOISSON}, we see if we choose $Z$ such that
\begin{equation*}
Z^{-1} \left(\frac{1}{\mf{L}_y^2} + \left|\exp\left(-\mf{L}_y^2\left(1-\cos\left(\frac{v}{\mf{L}_y}\right)\right)\right) - e^{-v^2}\right|\right) = o(\mf{L}_y^{-1})
\end{equation*}
then by Proposition \ref{CONTRA}, ii) in Proposition \ref{MAINPROP} cannot hold, and thus for any $T \geq 1$ we have
\begin{equation*}
|\phi_x(v,-v) - \chi(v,-v)| \ll Z\left(\int_{[T^{-3},T]^2} \frac{|\phi_{x,y}(t,t')-\chi(t,t')|}{|tt'|}dtdt' + T^{-1}\right).
\end{equation*}
%Now, suppose $\left|\frac{v}{\mf{L}}\right| \leq \log_2^{-\frac{1}{4}} x$. Then 
%\begin{align*}
%&\left|\exp\left(-\mf{L}_y^2\left(1-\cos\left(\frac{v}{\mf{L}_y}\right)\right)\right) - e^{-\frac{1}{2}v^2}\right| = e^{-\frac{1}{2}v^2}\left|\exp\left(\frac{1}{2}v^2-\mf{L}_y^2\left(1-\cos\left(\frac{v}{\mf{L}_y}\right)\right)\right)-1\right| \\
%&\ll e^{-\frac{1}{2}v^2} \mf{L}_y^2\left|\frac{v}{\mf{L}_y}\right|^4 \ll e^{-\frac{1}{2}v^2}v^4\mf{L}_y^{-2}.
%\end{align*}
%Conversely, when $\left|\frac{v}{\mf{L}_y}\right| > \log_2^{-\frac{1}{4}} x$ we have
%\begin{align*}
%&\left|\exp\left(-\mf{L}_y^2\left(1-\cos\left(\frac{v}{\mf{L}_y}\right)\right)\right) - e^{-v^2}\right| \leq \exp\left(-\mf{L}_y^2\left(1-\cos\left(\frac{v}{\mf{L}_y}\right)\right)\right) + e^{-v^2} \\
%&\ll \exp\left(-\frac{2}{3}\mf{L}_y^2\left|\frac{v}{\mf{L}_y}\right|\right) + e^{-v^2} \ll \exp\left(-\mf{L}_y^{\frac{3}{2}}\right).
%\end{align*}
(Here we have used the fact that $e^{-\frac{1}{21}\beta \log \beta} \ll \mf{L}^{-2}$.) Inserting the estimates for the term in absolute value derived in the proof of Proposition \ref{PHIEST} and noting that $\mf{L}_y \asymp \mf{L}_x$, we see that it suffices to choose $Z$ such that
\begin{equation*}
Z \asymp \begin{cases} (\log_4 x)\mf{L}^{-1} & \text{ if $\left|v\right| \leq \log_2^{\frac{1}{4}}x$} \\ (\log_4 x)\left(\mf{L}^{-2} + \exp\left(-\mf{L}_y^{\frac{3}{2}}\right)\right) &\text{ otherwise}. \end{cases} 
\end{equation*}
As such, according to the argument in the proof of Theorem \ref{EKTHM}, we get
\begin{equation*}
|\phi_{x,y}(v,-v) - \chi(v,-v)| \ll \begin{cases} (\log_4 x)(\log_3 x)(\log_2 x)^{-\frac{1}{4}} \mf{L}^{-1} &\text{ if $\left|v\right| \leq \log_2^{\frac{1}{4}}x$} \\ (\log_4 x)(\log_3 x)(\log_2 x)^{-\frac{1}{4}}\left(\mf{L}^{-2} + \exp\left(-\mf{L}_y^{\frac{3}{2}}\right)\mf{L}\right) &\text{ otherwise}. \end{cases}
\end{equation*}
Inserting these estimates into our expression for $\mc{E}$ gives
\begin{align*}
\mc{E} &\ll \mf{L}^{-1}\left(\int_{\left|v\right| \leq \mf{L}^{\frac{1}{2}}} + \int_{\left|v\right| > \mf{L}^{\frac{1}{2}}} \right) |\phi_{x,y}(v,-v)-\chi(v,-v)|dv \\
&\ll (\log_4 x)\left(\mf{L}^{-1}(\log_3 x)(\log_2 x)^{-\frac{1}{4}} + \mf{L}^2 \exp\left(-\mf{L}_y^{\frac{3}{2}}\right)\right) \\
&\ll \mf{L}^{-1}(\log_4 x)(\log_3 x)(\log_2 x)^{-\frac{1}{4}}.
\end{align*}
Since $e^{-\frac{j^2}{4\mf{L}_y^2}} \gg (\log_4 x)(\log_3 x)(\log_2 x)^{-\frac{1}{4}} $ when $|j| \leq \left(1-\e\right)\sqrt{\log_2 y\log_3 x}$, we indeed have $\mc{E} = o\left(e^{-\frac{j^2}{4\mf{L}_y^2}} \mf{L}^{-1}\right)$ in this case. Since $\log_2 x = (1+o(1)) \log_2 y$, this completes the proof. 
\end{proof}
\section*{Acknowledgments}
The author would like to thank his Ph.D supervisor John Friedlander, as well as Asif Zaman and Oleksiy Klurman for their encouragement and interest in the results of this paper. He is indebted to Youness Lamzouri for his helpful comments regarding how to improve the paper's exposition. Finally, the author thanks Maksym Radziwi\l \l $\ $ and Sary Drappeau for providing additional references.
\bibliographystyle{plain}
\bibliography{chowla}
\end{document}